\documentclass[12pt,a4paper]{article}
\usepackage[usenames,dvipsnames]{xcolor}
\usepackage{epsfig}
\usepackage{epsfig,amssymb,color, graphics}

\def\II{{\mathbb I}}   
   
  \def\SS{{\mathbb S}} 
   
 \def\ZZ{{\mathbb Z}}
\newtheorem{theo}{Theorem}
\newtheorem{lemm}{Lemma}[section]
\newtheorem{defi}[lemm]{Definition} 
\newtheorem{collary}[lemm]{Corollary} 
\newtheorem{remark}[lemm]{Remark} 

\newtheorem{prop}[lemm]{Proposition}
\newenvironment{demo}{{\bf Proof: }}{${}$\hfill $\diamond$ \medskip}

\def\al{\alpha}

\def\Om{\Omega}
\def\ga{\gamma}    
\def\Ga{\Gamma}
\def\de{\delta}
\def\De{\Delta}
\def\vp{\varphi}
\def\la{\lambda}

\def\si{\sigma}
\def\Si{\Sigma}
\def\ep{\varepsilon} 
\def\nd{\noindent}
\def\R{\mathbb R}
\def\D{\mathbb D}
\def\SS{\mathbb S}
\def\A{\mathbb A}
\def\Z{\mathbb Z}

\begin{document}
\title{Topological classification of Morse-Smale
diffeomorphisms without heteroclinic curves on $3$-manifolds}
\author{Ch.~Bonatti\and V.~Grines\and F. Laudenbach\and O.~Pochinka
\thanks{The dynamics and necessary and sufficient conditions of the topological conjugation of Morse-Smale
diffeomorphisms without heteroclinic curves on $3$-manifolds (Sections 2, 4) were investigated with the support of the Russian Science Foundation  (project 17-11-01041). The construction of the compatible foliations and topological research (Sections 3, 5) partially supported by RFBR (project nos. 15-01-03687-a, 16-51-10005-Ko\_a), BRP at the HSE (project 90) in 2017 and by ERC Geodycon.}}

\date{ }
\maketitle 
\sloppy \tableofcontents
\begin{abstract} 
{We show that, up to topological conjugation,} the equivalence class of a Morse-Smale
diffeomorphism without heteroclinic curves on a $3$-manifold is completely defined by an embedding of two-dimensional stable and unstable heteroclinic laminations to a characteristic space.
\end{abstract}
\section{Introduction and formulation of the result}\label{I}
In 1937 A. Andronov and L. Pontryagin \cite{AP} introduced the notion of a {\it rough system} {of} differential {equation} given in a bounded part of the plane, that is {a system which preserves} its qualitative properties under parameters variation if the variation is small enough. They proved that {the flows generated by such systems are exactly the flows having the} following properties:
\begin{enumerate}
\item the set of fixed points and periodic orbits is finite and all its elements are hyperbolic;
\item there are no separatrices going from {one saddle to itself or to another one}; 
\item all $\omega$- and $\alpha$-limit sets are contained in the union of fixed points and periodic orbits (limit cycles).
\end{enumerate}
The above description characterizes the rough flows on the two-dimensional sphere also. After  A. Mayer \cite{Mai} 
in 1939,  a similar result 
holds true on the 2-torus for 
flows  having a  closed section and no equilibrium states. 
 A. Andronov and L. Pontryagin  have shown also in \cite{AP} that the set of  the 
rough flows is  dense in the space of $C^{1}$-flows\footnote{This statement  was not explicitly formulated  in 
\cite{AP} and  was  mentioned for the first time in papers by E. Leontovich \cite{Le} and M. Peixoto \cite{Pe0}. G. 
Baggis \cite{Ba} in 1955 made explicit some details of the proofs which were not completed in \cite{AP}.}. 
In 1962 M. Peixoto proved (\cite{Pe1}, \cite{Pe2}) that the properties 1-3 are necessary and sufficient for a flow 
on any orientable surface of genus greater than zero to be structurally stable. He proved the density for these flows
 as well.   Direct generalization of the properties of rough flows on surfaces leads to the following class of dynamical
  systems continuous or discrete, that is, flows or diffeomorphisms (cascades).

\begin{defi} \label{oM-S} A smooth dynamical system  given on an  $n$-dimensional manifold  $(n\geq 1)$ 
$M^n$ is called  Morse-Smale if:
\begin{enumerate} 
\item its non-wandering set consists of {a finite number of} fixed points and periodic orbits {where} each of 
them is hyperbolic;
\item the stable and unstable manifolds  $W^s_p$, $W^u_q$ of any pair of  non-wandering points $p$ and $q$ 
intersect transversely. 
\end{enumerate}
\end{defi}

Let $M$ be a given closed 3-dimensional manifold and $f:M\to M$ be a Morse-Smale diffeomorphism. 

For $q=0,1,2,3$ denote by $\Omega_q$ the set of all periodic points of $f$ with $q$-dimensional  unstable 
manifold. Let $\Om_f$ be the union of all periodic points. Let us represent the dynamics of $f$ in the form 
``source-sink'' {in the} following way. Set  
$$A_f=W^u_{\Omega_0\cup\Omega_1},\,\,R_f=W^s_{\Omega_2\cup\Omega_3},\,\,V_f=M
\smallsetminus(A_f\cup R_f).$$ 

We recall that a compact set $A\subset M$ is said to be an {\it attractor} of $f$ if there is a compact 
neighborhood $N$ of the set $A$ such that $f (N) \subset int~N$ and $A=\bigcap\limits_{n\in\mathbb{N}}f^n(N)$; and $R\subset M$
is said to be a {\it repeller} of $f$ if it is an attractor of $f^{-1} $.

By \cite[Theorem 1.1]{GrMePoZh}  the set {$A_f$ (resp. $R_f$)  is  an attractor (resp.  a repeller)}
of $f$ whose topological dimension is equal to 0 or 1. By \cite[Theorem 1.2]{GrMePoZh} the set $V_f$ is a connected 3-manifold and $V_f=W^s_{A_f\cap \Omega_f}\smallsetminus A_f=W^u_{R_f\cap \Omega_f}\smallsetminus R_f$. Moreover, the quotient $\hat V_f=V_f/f$ is a closed connected 3-manifold and when $\hat V_f$ is orientable, 
then it is either irreducible  or diffeomorphic to $\mathbb S^2\times\mathbb S^1$. The natural projection $p_{_f}:V_f\to\hat V_f$  {is an infinite cyclic covering. Therefore, there is a natural epimorphism from the the first homology group of $\hat V_f$ to $\ZZ$, $$\eta_{_f}:H_1(\hat V_f;\ZZ)\to\ZZ,$$ defined as follows: if $\ga$ is a path in $V_f $  joining $x$ to $f^n(x)$, $n\in \ZZ$, then $\eta_{_f}$ maps the homology class of the cycle $p_{_f}\circ\ga$ to $n$.}

The intersection with $V_f$ of the $2$-dimensional stable manifolds of the saddle points of $f$ is an invariant $2$-dimensional lamination $\Ga^s_{f}$, with finitely many leaves, and which is closed in $V_f$. Each leaf of this lamination is obtained by removing from a stable manifold its set of intersection points with the $1$-dimensional unstable manifold; this intersection is at most countable. As $\Ga^s_{f}$ is invariant under $f$, it 
descends 
to the quotient in a compact 
$2$-dimensional lamination $\hat\Ga^s_{f}$ on $\hat V_f$. Note that each $2$-dimensional stable manifold is a plane on which $f$ acts as a contraction, so that the quotient by $f$ of the punctured  stable manifold is either a torus or a Klein bottle.  Thus the leaves of $\hat\Ga^s_{f}$ are either tori or Klein bottles which are punctured along at most countable set.

One defines in the same way the unstable lamination $\hat \Ga^u_{f}$ as the quotient by $f$ of the intersection with $V_f$ of the $2$-dimensional unstable manifolds. The laminations $\hat\Ga^s_{f}$ and $\hat \Ga^u_{f}$ are transverse.

\begin{defi} The sets $\hat{\Gamma}^s_f$ and  $\hat{\Gamma}^u_f$ are called {the} two-dimensional stable and unstable laminations associated  with the diffeomorphism $f$.
\end{defi}
A precise definition of what a lamination is  will be given in Definition \ref{lamination}. 
\begin{defi} \label{s-s-s} The collection $S_{f}=(\hat V_{f},\eta_{_{f}},\hat{\Gamma}^s_{f},\hat{\Gamma}^u_{f})$ is called the scheme of the diffeomorphism $f$.
\end{defi}
\begin{defi} \label{eqg03} {The schemes} $S_f$ and  $S_{f'}$ of  {two} Morse-Smale diffeomorphisms 
$f,f':M\to M$ are {said to be} equivalent if there is a homeomorphism  $\hat\varphi:\hat V_f\to\hat V_{f'}$ 
with following properties:\\
{\rm (1)}\quad $\eta_{_f}=\eta_{_{f'}}\hat\varphi_*$;\\
{\rm (2)}\quad $\hat\varphi(\hat{\Gamma}^s_{f})=\hat{\Gamma}^s_{f'}$ and 
 $\hat\varphi(\hat{\Gamma}^u_{f})=\hat{\Gamma}^u_{f'}$, meaning that $\hat\vp$ maps leaf to leaf.
\end{defi}  

Using the above notion of a scheme in a series of papers by Ch. Bonatti, V. Grines, V. Medvedev, E. Pecou, O. Pochinka \cite{BoGr}, \cite{BoGrMePe3}, \cite{BoGrPo2005}, \cite{BoGrPo2}, the problem of classification
up to topological conjugacy of Morse-Smale diffeomorphisms on 3-manifolds has been solved in some particular cases. Recall that two diffeomorphisms $f$ and $f^{\prime}$ of $M$  are said to be {\it topologically conjugate} if there is a homeomorphism $h:M\to M$ which satisfies $f'h=hf$.   

In the present {article,} we give the topological classification of the  Morse-Smale diffeomorphisms 
belonging to the subset $G(M)$ of the Morse-Smale diffeomorphisms $f:M\to M$ which have no heteroclinic
 curves (see Section \ref{section2}). According to \cite{BGMP},  when  the ambient manifold is orientable, 
then it is either sphere $\mathbb S^3$ or the connected sum of a finite number copies of 
$\mathbb S^2\times\mathbb S^1$. 
\begin{theo} \label{t.invariant} {Two Morse-Smale diffeomorphisms in $G(M)$ are topologically conjugate} 
if and only if their schemes are equivalent.
\end{theo}

The structure of the paper is the following:
\begin{itemize}
\item In Section \ref{section2} we describe the dynamics of Morse-Smale diffeomorphisms and their space of wandering orbits. 
\item In Section \ref{section-compatible} we construct a compatible system of neighborhoods, which is a key point for the  construction of a conjugating homeomorphism. 
\item In Section \ref{IV} we construct a conjugating homeomorphism.     
\item  Section \ref{V} is an appendix of 3-dimensional topology. We prove there some topological lemmas 
which are used in Section 4. 
\end{itemize}

\section{Dynamics of diffeomorphisms {in the class} $G(M)$}\label{section2}
In this section we introduce some notions connected with Morse-Smale diffeomorphisms on 3-manifold $M$.  
More detailed information on Morse-Smale diffeomorphisms is contained in \cite{grin} for example.

\begin{figure}[h]
\centerline{\includegraphics[width=11cm,height=6cm]{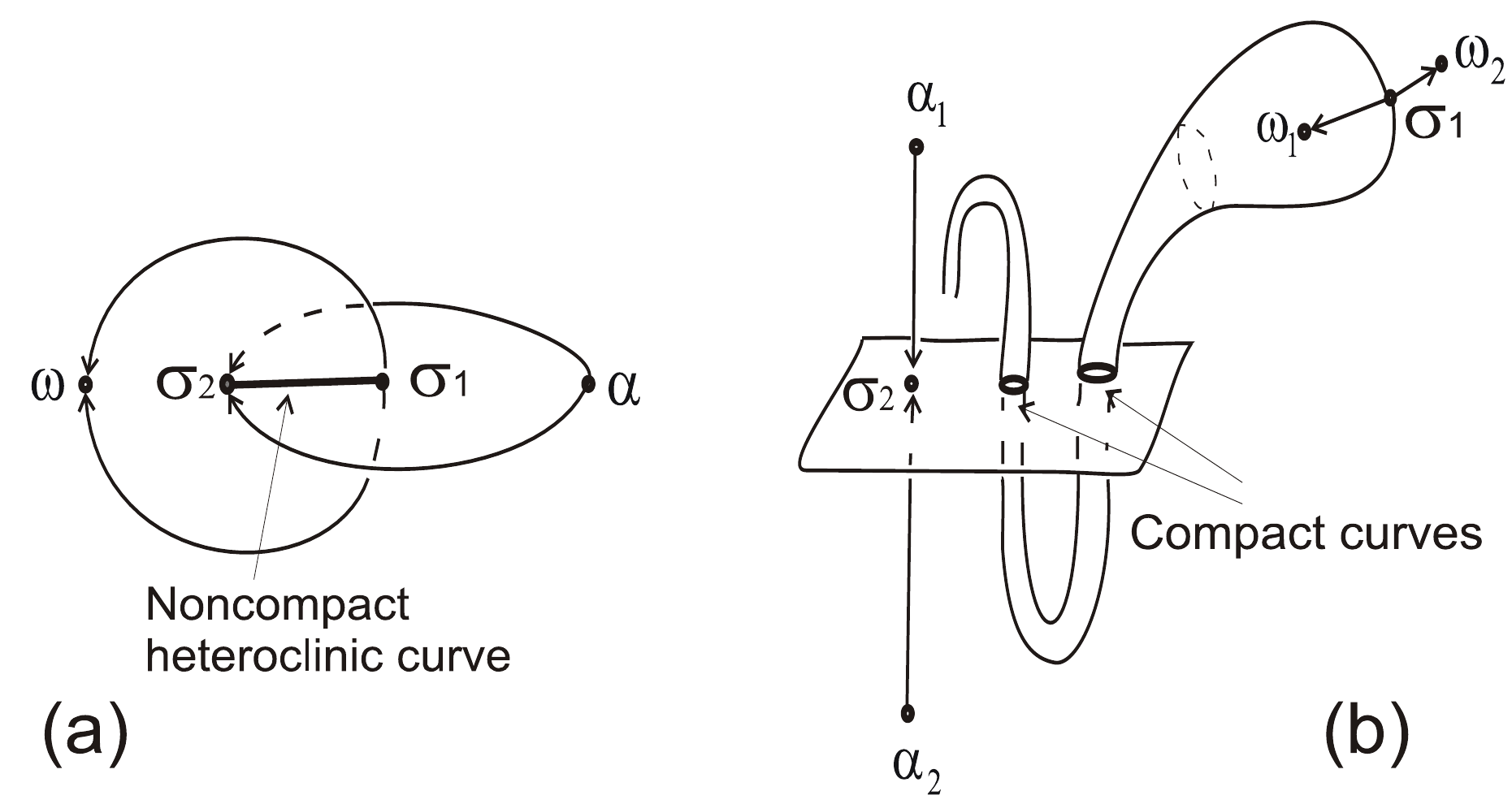}}\caption{Heteroclinic curves} \label{pic3}
\end{figure}

Let $f:M\to M$ be a Morse-Smale diffeomorphism. If $x$ is a periodic point its {\it Morse index} is the 
dimension of its unstable manifold $W^u_x$; the point $x$ is called a {\it saddle point} when its two invariant 
manifolds have positive dimension, that is, its Morse index is not extremal. A {\it sink point} has Morse index
 $0$ and a {\it source point} has Morse index {$3$}. The following notions are key concepts for describing 
 how the stable manifolds of saddle points intersect the unstable ones.
If $x, y$ are distinct 
saddle points of $f$ and $W^u_x\cap W^s_y\neq\emptyset$, then:
\begin{itemize}
\item[-] {if $dim~W^s_x < dim~W^s_y$, any connected component of $W^u_x\cap W^s_y$ is 1-dimensional and called a {\it heteroclinic curve}} (see figure \ref{pic3});
\item[-] {if $dim~W^s_x=dim~W^s_y$, the set $W^u_x\cap W^s_y $ is  countable;  each  of its points is called a {\it heteroclinic point}}; the orbit of a heteroclinic point is called a {\it heteroclinic orbit}.
\end{itemize}
\begin{figure}[h]
\centerline{\includegraphics[width=9cm,height=6cm]{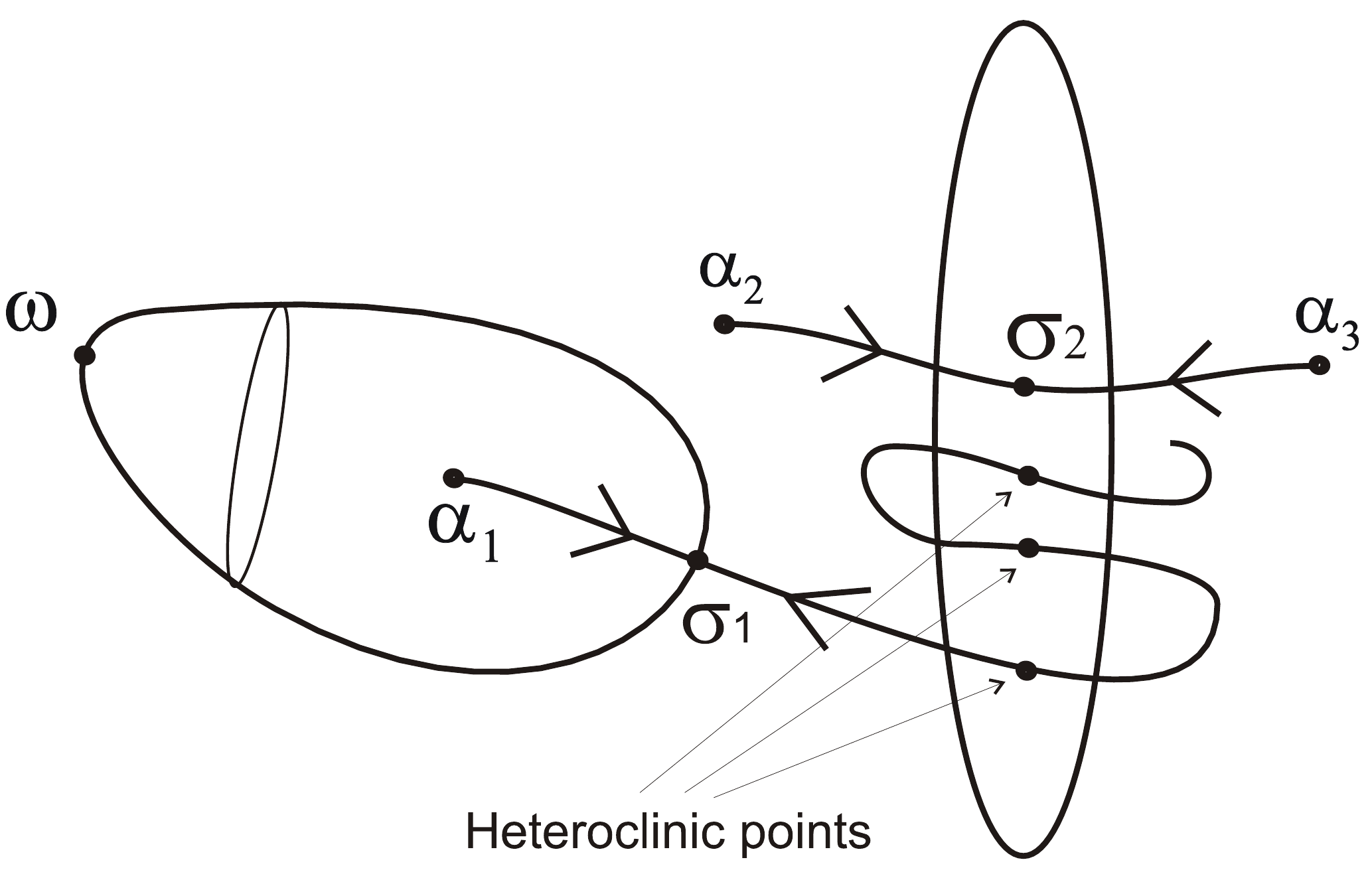}}\caption{Heteroclinic points} \label{pic4}
\end{figure}
According to S. Smale \cite{S3}, it is possible to define a partial order in the set of saddle points of a given Morse-Smale diffeomorphism $f$ as follows: for different periodic orbits $p\neq q$,  one sets $p\prec q $ if  and only if  $W^u_q \cap W^s_p \neq\emptyset${.}  Smale proved that this relation is a partial order. In that case, it follows from \cite[Lemma 1.5]{Pa} that there is a sequence of different periodic orbits $p_0, \dots, p_n $ satisfying the following conditions: $p_0=p $, $p_n=q $ and $p_i\prec p _ {i+1} $.  The sequence $p_0, \dots, p_n $ is said to be an {\it $n$-chain connecting $p$ to $q$}. The length of the longest chain connecting $p$ {to} $q$ is denoted by $beh(q|p)$. {If  $W^u_q \cap W^s_p =\emptyset$, we pose $beh(q|p)=0$}. For a subset $P$ of the periodic orbits let us set $beh(q|P)=\max\limits_{p\in P}\{beh(q|p)\}$. The present paper is devoted to studying Morse-Smale diffeomorphisms  in dimension 3 which have no heteroclinic curves. We recall from the introduction that this class of diffeomorphisms is denoted by $G(M)$. 
Let $f\in G(M)$. It follows from \cite{GrMePoZh} that if the set $\Omega_2$ is empty then 
$R_f$ consists of a unique source. If $\Omega_2\neq\emptyset$,   denote by $n$ the length of the longest 
chain connecting {two points of}  $\Omega_2$. Divide the  set $\Omega_2$ {into} $f$-invariant
 disjoint parts $\Si_0,\Si_1,\dots,\Si_n$ using the rule: $beh(q|(\Omega_2\smallsetminus q))=0$ for 
 each orbit $q\in \Si_{0}$ and $beh(q|\Si_i)=1$ for each orbit $q\in \Si_{i+1},~i\in\{0,\dots,n-1\}$. 
{Since $\Omega_1$ for $f$ is $\Omega_2$ for $f^{-1}$,} then it is possible to divide the periodic orbits of the set $\Omega_1$ into parts in a  similar way. The absence of heteroclinic curves means that there are no chains connecting a saddle from $\Omega_2$ with a saddle from $\Omega_1$. Thus we explain all material for $\Omega_2$ and say that all is similar for $\Omega_1$. 
 
Set $W^u_i:=W^u_{\Si_i}$, $W^s_i:=W^s_{\Si_i}$. Then, $R_f:=\bigcup
\limits_{i=0}^n cl(W^s_i)$, where $cl(\cdot)$ stands for the closure of~$(\cdot)$. We now specify what 
 a lamination {is} and which sort of regularity it may have. There 
are different possible notations. Here, we use the one which is given in \cite[{Definition 1.1.22}]{candel}.
 
\begin{defi} \label{lamination} Let $X$ be a $n$-dimensional and $Y\subset X$ be a closed subset. Let $q$ be an integer $0<q< n$.  A codimension-$q$ \emph{lamination} with support $Y$ is a decomposition $Y= \bigcup\limits_{j\in J}L_j$ into pairwise disjoint smooth $(n-q)$-dimensional connected manifolds $L_j$, which are called the leaves. The family $L=\{L_j, j\in J\}$ is said to be a $C^{1,0}$-lamination\footnote{There 
are different possible notations. Here, we use the one which is given in \cite[{Definition 1.1.22}]{candel}.}
 if for every point $x\in Y$ the following conditions hold:
 \begin{itemize}
 \item[\rm (1)] There are an open neighborhood $U_x\subset X$ of $x$ and a homeomorphism 
$\psi:U_x\to\mathbb R^n$ such that $\psi$ maps every \emph{plaque}, that is a connected component of 
 $U_x\cap L_j$,  into a codimension-$q$ subspace 
  $\{(x^1,\dots,x^n)\in\mathbb R^n \mid x^{n-q+1}=c^{n-q+1},\dots,x^{n}=c^{n}\}$. If $Y=X$ 
  one says that $\mathcal{L}$ is foliation.
\item[\rm (2)] 
The tangent plane field   $TY:= \bigcup\limits_{j\in J}TL_j$ exists on $Y$ and is continuous.
\end{itemize}
\end{defi}

By definition,  two points belong to the same {\it leaf} of a lamination if they are linked by a path which is covered by finitely many plaques.

By abuse, a lamination and its support are generally denoted in the same way. We recall the $\lambda$-Lemma  in the strong form which is  proved in \cite[Remarks p. 85]{PaMe1998}. 

\begin{lemm}{\rm (\bf $\lambda$-lemma.)} \label{lambda} Let $f:X\to X$ be a diffeomorphism of 
 an $n$-manifold, and let $p$ be a fixed point of $f$.
Denote $W^u_{p}$ and $W^s_{p}$ be  the unstable and stable manifold respectively; 
say $\dim~W^u_{p}=m$, $0<m<n$. Let $B^s$ be a compact subset of $W^s_p$ (containing $p$ or not) and
 let $F: B^s\to C^1(\mathbb D^m, X) $ be a continuous family of embedded closed $m$-disks of class
  $C^1$ 
  transverse to $W^s_p$ and meeting $B^s$;
   set $F(x):= D^u_x$. Let $D^u\subset W^u_p$ be a compact $m$-disk 
and let $V\subset X$ be a compact $n$-ball such that $D^u$ is a connected component of $W^u_p\cap V$.  Then, when $k$ goes to $+\infty$, the sequence $f^k(D^u_x)\cap V$ converges to $D^u$ in the $C^1$ topology uniformly for $x\in B^s$.
\end{lemm}

Notice that it is important for applications that $B^s$ may not contain the point $p$.
\begin{figure}[h]
\centerline{\includegraphics[width=11cm,height=8cm]{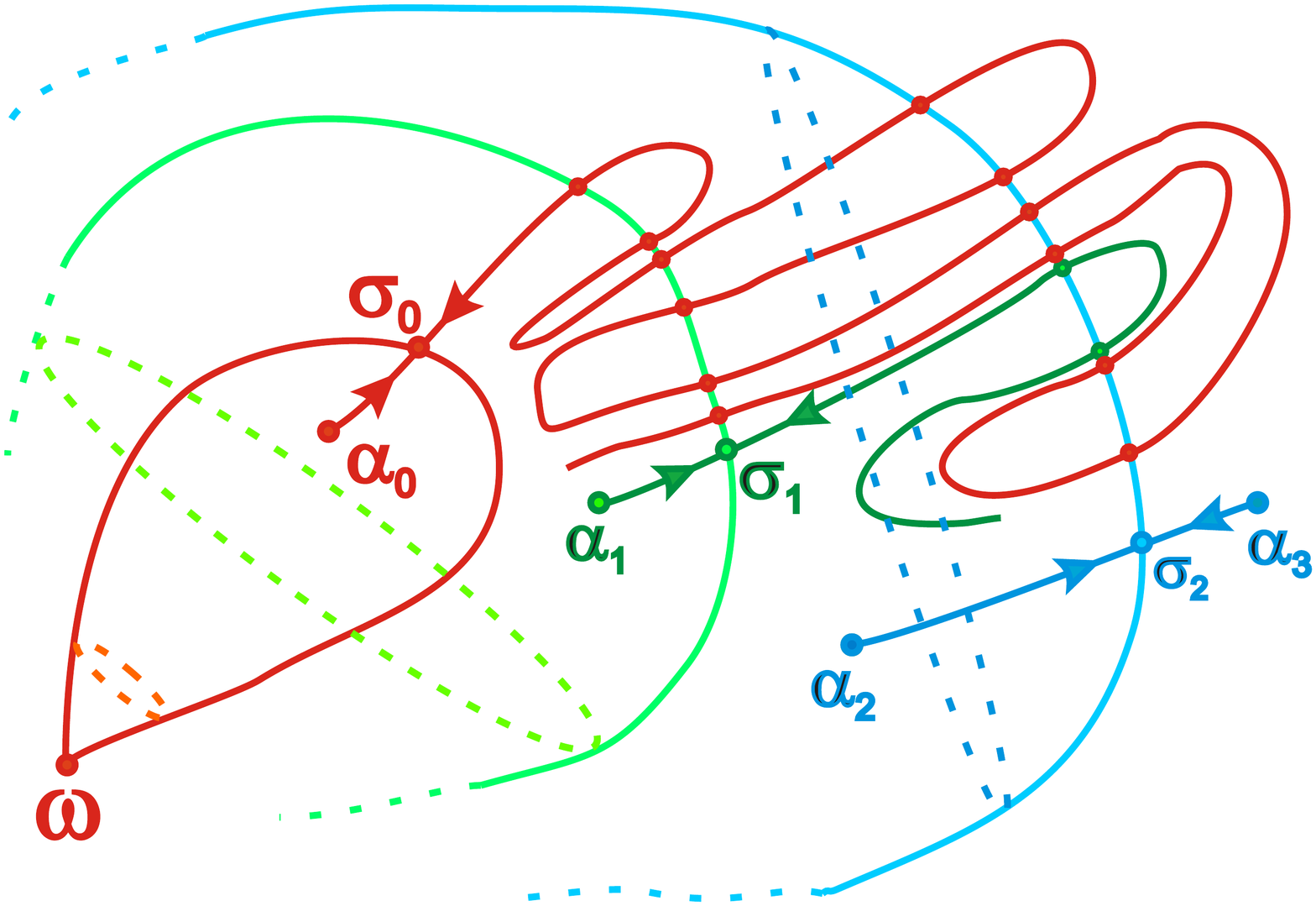}}\caption{A phase portrait of a diffeomorphism from the class $G(M)$.} \label{faz}
\end{figure}
Going back to our setting, a first application of the $\la$-lemma is that we have $W^u_i\subset cl(W^u_{i+1})$ and the closure of $cl(W^u_0) = W^u_0\cup \Om_0$. Moreover, $cl(W^u_n)\cap (M\smallsetminus \Om_0)$ is a  $C^{1,0}$-lamination of codimension one. From this one derives that $\hat L^u_f$ is also a $C^{1,0}$-lamination. Here is a typical example.

\begin{figure}[h]
\centerline{\includegraphics[width=17cm,height=6cm]{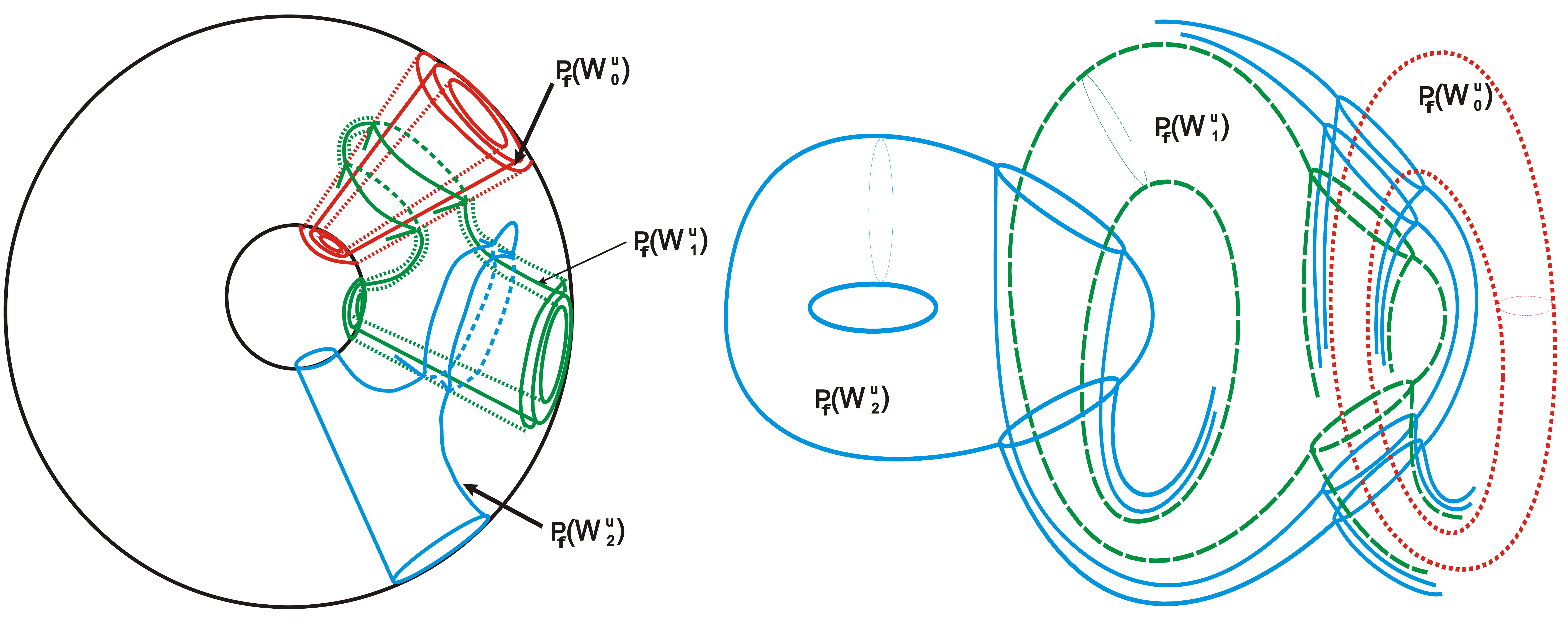}}
\begin{center}
\caption{A lamination associated with the diffeomorphism $f\in G(M)$ whose phase portrait is pictured in  Figure \ref{faz}.} \label{sch}
\end{center}
\end{figure}

On Figure \ref{faz} there is a phase portrait of a diffeomorphism $f\in G(M)$ whose non-wandering set 
$\Omega_f$ consists of fixed points: one sink $\omega$, 
three saddle points $\Si_0=\sigma_0,\ \Si_1=\sigma_1,\ \Si_2=\sigma_2$
with two-dimensional unstable manifolds and four sources $\alpha_0,\alpha_1,\alpha_2,\alpha_3$. 
We will illustrate all further proofs with this diffeomorphism. For this case $V_f:=W^s_\omega\smallsetminus\{\omega\}$. As the restriction of $f$ to the basin $W^s_\omega$ of $\omega$ is topologically conjugate to any  homothety,  $\hat V_f$ is diffeomorphic to $\mathbb S^2\times\mathbb S^1$. As $f|_{W^u_i}$ is topologically conjugate  to a homothety
 then  $(W^u_i\smallsetminus\Si_i)/f$ is diffeomorphic to the 2-torus; but this torus does not embed  to $\hat V_f$, except when $i=0$. On Figure \ref{sch} there is {the lamination associated with} the diffeomorphism $f\in G(M)$ whose phase portrait is {on Figure 3}. {On the left, the lamination is embedded in} $\mathbb S^2\times\mathbb S^1$
 which is seeen as the double of $\mathbb S^2\times \mathbb D^1$. 

We are going to show that the topological classification of diffeomorphisms in the class $G(M)$ reduces to classifying some appropriate  laminations {$\hat \Ga^u_f$ and $\hat \Ga^s_f$}. The technical key to the proof consists of constructing special foliations in some neighborhoods of the laminations.   

\section{Compatible foliations}\label{section-compatible}
Let $f\in G(M)$. Recall that we divided the set $\Omega_2$ into the $f$-invariant parts $\Si_0,\dots,\Si_n$. 
Using this partition, we explain how to construct compatible foliations (see Definition \ref{2dopsystem}) around $W^u_{\Omega_2}\cup W^s_{\Omega_2}$. {Similarly,} it is possible to construct compatible 
foliations around $W^s_{\Omega_1}\cup W^u_{\Omega_1}$. In what follows, we give ourselves  
four models of concrete hyperbolic linear isomorphisms $\mathcal E_{\kappa,\nu}  \in GL(\R^3),\,\kappa,\nu\in\{-,+\}$ given by the following formula: $$\mathcal E_{\kappa,\nu}(x_1,x_2,x_3)=(\kappa  2x_1,2x_2,\nu \frac{x_{3}}{4}).$$
The origin $O$ is the unique fixed point which is a saddle point with unstable manifold $W^u_O=Ox_1x_2$ and stable manifold $W^s_O=Ox_3$. If $\kappa=+$ (resp. $-$), the orientation of the unstable manifold
is preserved (resp. reversed), and similarly for the orientation of the stable manifold with respect to $\nu$. 
 We refer to each of them as the {\it canonical diffeomorphism}; it will be  denoted by $\mathcal E$  ignoring the sign.
 For $p\in \Omega_2$, let $per(p)$ denote the period of {$f$ at {$p$.} 
 
\begin{defi} \label{adop} A neighborhood $N_p$ of a saddle point $p\in\Omega_2$ is called \emph{linearizable} if there is a homeomorphism $\mu_p:N_p\to {\mathcal N}$ which conjugates the diffeomorphism  $f^{per(p)}\vert_{{N}_p}$ to the  canonical diffeomorphism $\mathcal E\vert_{{\mathcal N}}$.
\end{defi}
According to the local topological {classification of hyperbolic fixed point}  \cite[Theorem 5.5]{PaMe1998}, {every $p\in \Om_2$ has a linearizable neighborhood $N_p$.} For $t\in (0,1)$, set $\mathcal N^t:=\{(x_1,x_2,x_3)\in\mathbb{R}^3\mid 
-t< (x_1^2+x_2^2)x_{3}<  t\}$  and $\mathcal N:=\mathcal N^1$. 
The set   $\mathcal N^t$ is invariant by the canonical diffeomorphism $\mathcal E$.
By \cite{S3}, $W^s_p$ and $W^u_p$ are smooth submanifolds of $M$. 
The {\it boundary} of $\mathcal N$ is the surface in $\mathbb R^3$ defined by the equations $(x_1^2+x_2^2)x_3=\pm 1$. The open manifold $N_p$ has a similar boundary in $M$ denoted by $\partial N_p$.  This boundary is formed by points which are not in $N_p$ but are limit points of arcs in $N_p$\,; it is distinct from its closure as a subset of $M$. Clearly, {the linearizing homeomorphism} $\mu_p$ extends to $\partial N_p$.
For each $i\in \{0,\ldots, n\}$, choose some $p\in\Si_i$ and $\mu_p$ conjugating $f^{per(p)}$ to $\mathcal E\vert_{{\mathcal N}}$. Then, for $k\in\{1,\dots,per(p)-1\}$ define $\mu_{f^k(p)}$ 
so that the next formula holds for every $x\in N_{f^{k-1}(p)}$:
$$\mu_{f^k(p)}(f(x))=\mu_{f^{k-1}(p)}(x).$$ 
{We define  a pair of transverse foliations 
$\left(\mathcal{F}^u,\mathcal{F}^s\right)$  in $\mathcal N$} in the following way: 

-- the leaves of $\mathcal{F}^u$ are the fibres in $\mathcal N$ of the projection  $(x_1,x_2,x_3)\mapsto x_3$;

-- the leaves of $\mathcal{F}^s$ are the fibres in $\mathcal N$ of the projection $(x_1,x_2,x_3)\mapsto (x_1,x_2)$.

\nd By construction, $W^u_O$ and  $W^s_O$ are  leaves  of $\mathcal F^u$  and  $\mathcal F^s$ respectively. Let $N_i$ denote the union $\bigcup\limits_{p\in\Si_i}N_p$. This is an $f$-invariant neighborhood of $\Si_i$. Let $\mu_i:N_i\to \mathcal N$ be the map whose restriction to $N_p$ is $\mu_p$. Thus, taking the pullback of them by $\mu_i$ gives  a pair of $f$-invariant foliations $\left({F}^u_{i},{F}^s_{i}\right)$ on $N_i$ which are said to be {\it linearizable}. By construction, $W^u_i$ and $W^s_i$ are made of leaves of $F^u_i$ and $F^s_i$ respectively.  Sometimes we want to deform the linearizable neighborhood  $N_p$ by {\it shrinking}. Observe that the homotheties of ratio $\rho\in (0,1)$ act on $\mathcal N$ preserving 
$\mathcal F^u$ and $\mathcal F^s$ and map $\mathcal N$ to $\mathcal N^{\rho^3}$. By conjugation, similar contractions $c_\rho $  are available  in $N_p$ for every $p\in\Omega_2$. The neighborhood  $c_\rho(N_p)$ is said to be obtained from $N_p$ by {\it shrinking}.
 
\begin{lemm} \label{Ge} For every $\rho\in (0,1)$, the shrunk neighborhood $c_\rho (N_p)$ is linearizable. Generically, the boundary of $c_\rho(N_p)$ does not contain any heteroclinic point.
\end{lemm}
\begin{demo} For a given $\mu_p$, we define $\mu_p^\rho$ as follows: its domain is $\mu_p^{-1}(\mathcal N^{\rho^3})$ and, on this domain, it is defined by $\mu_p^\rho= c_\rho^{-1}\circ\mu_p$. Its range is 
$\mathcal N$. Since the heteroclinic points form a countable set, for almost every $\rho\in (0,1)$ the boundary of the domain of $ \mu_p^\rho$ avoids the heteroclinic points.
\end{demo}

Observe that the canonical diffeomorphism and the contraction $c_\rho$ keep both foliations invariant. Recall the $f$-invariant partition  $\Om_2=\Si_0\sqcup\Si_1\sqcup\ldots\sqcup\Si_n$.
Let us introduce the following notations:
\begin{itemize} 
\item[-] for any $t\in(0,1)$, set  $N^t_p:=\mu_p^{-1}(\mathcal N^t)$ and $N^t_i:=\bigcup\limits_{p\in\Si_i}N^t_p$\,;
\item[-] for any point $x\in N_i$, denote  ${F}^u_{i,x}$ (resp. ${F}^s_{i,x}$) the  leaf of the foliation  ${F}^u_{i}$ (resp. ${F}^s_{i}$) passing through $x$;
\item[-] for each point $x\in N_i$, set $x^u_i=W^u_i\cap F^s_{i,x}$ and $x^s_i=W^s_i\cap F^u_{i,x}$. 
Thus, we have $x=(x^u_i,x^s_i)$ in the coordinates defined by $\mu_i$. \end{itemize}
{We also introduce the {\it radial functions} $r^u_i, r^s_i: N_i\to [0,+\infty)$ defined by:
$$r^u_i(x)= \Vert \mu_i(x^u_i)\Vert^2\quad{\rm and} \quad r^s_i(x)= \vert \mu_i(x^s_i)\vert.$$
With this definition at hand, the neighborhood $N_i^t$ ot $\Si_i$ is defined by the inequality
$$r^u_i(x)r^s_i(x)<t\,.$$
Observe that the radial function $r^s_i$ endows each stable separatrix of $p\in \Si_p$ with a natural order which will be used  later in the proof of Theorem \ref{t.invariant}.}

\begin{defi}\label{2dopsystem} The linearizable neighborhoods $N_0,\dots,N_n$ are called \emph{compatible} if, for any $0\leq i<j\leq n$ and   $x\in N_{i}\cap N_{j}$, the following holds:
$$F^s_{j,x}\cap N_{i}\subset F^s_{i,x}\,\,\,\,and\,\,\,F^u_{i,x}\cap N_{j}\subset F^u_{j,x}.$$
\end{defi}
If linearizable neighborhoods are compatible, they remain so after some of them are shrunk. 
\begin{remark} {\rm The notion of compatible foliations is a modification of the admissible systems of tubular families introduced by 
J. Palis and S. Smale in \cite{Pa} and \cite{PS}.}
\end{remark}
We introduce the following notation:
\begin{itemize} 
\item[-]  For $i\in\{0,\dots,n\}$, set 
$A_{i}:=A_f\cup\bigcup\limits_{j=0}^i W_j^u,\hfill\break
V_i:= W^s_{A_i\cap\Om_f}\smallsetminus A_i,~~\hat V_i:= V_i/f$. Observe that $f$ acts freely 
on $V_i $ and denote the natural projection by $p_{_{i}}:V_{i}\to\hat V_{i}${.} 
\item[-] For $j,k\in\{0,\dots,n\}$ and $t\in(0,1)$, set $\hat W^s_{j,k}=p_k(W^{s}_{j}\cap V_k)$, $\hat W^u_{j,k}=p_k(W^{u}_{j}\cap V_k)$,  $\hat{N}^{t}_{j,k}=p_{k}(N^t_j\cap V_{k})${.}
\item[-]  $L^u:=\bigcup\limits_{i=0}^nW^u_i$, $L^s:=\bigcup\limits_{i=0}^nW^s_i$,
$L^u_i:=L^u\cap V_i$\,,
~$L^s_i:= L^s\cap V_i$\,,
$\hat{L}^u_i:=p_i(L^u_i)$, $\hat{L}^s_i:=p_i(L^s_i)$. 
\end{itemize} 
\begin{theo} \label{compatible} For each diffeomorphism $f\in G(M)$ there exist compatible linearizable neighborhoods of all saddle points whose Morse index is $2$. 
\end{theo}
\begin{demo} The proof consists of three steps. 

\nd {\bf Step 1.} Here, we prove the following claim.
\begin{lemm} There exist $f$-invariant neighborhoods $U^{s}_{0},\dots,U^{s}_{n}$ of the sets $\Si_{0},\dots,\Si_{n}$ respectively,  equipped with  two-dimensional $f$-invariant  foliations   ${F}^{u}_{0},\dots,{F}^{u}_{n}$ {of class $C^{1,0}$} such that  the following properties hold for each $i\in\{0,\dots,n\}$:
\begin{enumerate}
\item[\rm (i)] the unstable manifolds $W^u_i$ are leaves of the foliation ${F}^{u}_{i}$ and each leaf of the foliation ${F}^{u}_{i}$ is transverse to  $L^s_i$;
 \item[\rm (ii)] for any  $0\leq i <k\leq n$ and $x\in U^s_{i}\cap U^s_{k}$, we have the inclusion $F^u_{k,x}\cap {U}^s_{i}\subset F^u_{i,x}$. 
\end{enumerate}
\end{lemm}
\begin{demo} Let us prove this by a decreasing induction on $i$ from $i=n$ to $i=0$. For $i=n$, it follows from the 
definition of $V_n$ that $(W^s_{n}\smallsetminus\Si_{n})\subset V_{n}$. Since $f$ acts freely and properly on $W^s_n$, the quotient $\hat W^{s}_{n,n}$ is a smooth submanifold of $\hat V_{n}$; it  consists of finitely many 
circles. The lamination $\hat L^s_n$ accumulates on $\hat W^{s}_{n,n}$. Choose an open tubular neighborhood 
${\hat N}^{s}_{n}$ of  $\hat W^{s}_{n,n}$ in $\hat V_n$; denote its projection by $\pi^{u}_{n}:{\hat N}^{s}_{n}\to \hat W^{s}_{n,n}$. Its fibers form a 2-disc foliation $\{{d}^{u}_{n, x}\mid x\in \hat W^{s}_{n,n}\}$  {transverse to 
$\hat W^{s}_{n,n}$. Since  $\hat L^s_n$ is a $C^{1,0}$-lamination  containing $\hat W^{s}_{n,n}$, each plaque of 
$\hat W^{s}_{n,n}$ is the $C^1$-limit of any sequence  of  plaques approaching it $C^0$. Therefore, if the tube ${\hat N}^{s}_{n}$ is small enough, its fibers are  transverse  to $\hat L^s_n$.}
 
Set $U^{s}_{n}:=p_{n}^{-1}({\hat N}^{s}_{n})\cup W^u_{n}$. This is an open set of $M$ which carries a foliation $F^u_n$ defined by {taking the preimage of} the fibers of $\pi^u_n $ and by adding $W^u_{n}$ as extra leaves. This is the requested foliation satisfying (i) and (ii) for $i=n$. Notice that the plaques  of $F^u_n$ are smooth {and by the $\lambda$-lemma, for any compact disc  $B$ in $W^u_n$ there is $\varepsilon >0$ such that  every plaque  of $F^u_n$  which is $\ep$-close to $B$ in topology $C^0$ is also $\varepsilon $-close to $B$ in topology $C^1$. Hence, $F^u_n$ is a $C^{1,0}$-foliation.} 

For the induction, we assume the construction is done for every $j>i$ and we have to construct an $f$-invariant neighborhood $U^{s}_{i}$ of the saddle point{s in $\Si_{i}$} carrying an $f$-invariant 
foliation ${F}^{u}_{i}$ satisfying (i) and (ii). Moreover, by genericity the boundary  $\partial U^s_{j}$, $j>i$, is assumed to avoid all heteroclinic points. For  $j>i$, let $\hat{U}^{s}_{j,i}:=p_{i}(U^{s}_{j}\cap V_{i})$ and 
$\hat{{F}}^{u}_{j,i}:=p_{_{i}}({F}^{u}_{j}\cap V_{i})$. For the same reason as in the case $i=n$, the set $\hat W^{s}_{i,i}$ is a smooth submanifold of $\hat V_{i}$ consisting of circles. Choose a tubular neighborhood  $\hat N^{s}_{i}$ of $\hat W^{s}_{i,i}$ with a projection $\pi^{u}_{i}:{\hat N}^{s}_{i}\to \hat W^{s}_{i,i}$ whose fibers are 2-discs. Similarly, $(W^u_{i+1}\setminus\Si_i)\subset V_i$ and, hence, $\hat W^u_{i+1, i}$ is a compact submanifold, consisting of finitely many  tori or Klein bottles. The set $\hat L^u_i$ is  a compact lamination  and its intersection with $\hat W^{s}_{i,i}$ consists of a countable set of points which are the  projections of the heteroclinic points belonging to the stable manifolds  $W^{s}_{i}$. Actually, there is a hierarchy in $\hat L^u_i\cap \hat W^{s}_{i,i}$  which we are going to describe in more details.
 
Set $H_{k}:= \hat W^u_{i+k,i}\cap \hat W^s_{i,i}$ for $k>0$. Since $\hat W^u_{i+1, i}$ is compact, $H_1$ is a finite set: $H_1= \{h^1_1, ... , h^1_{t(1)}\}$.  We are given neighborhoods,  called {\it  boxes},  $B^1_\ell$, $\ell=1,..., t(1)$, about these points,   namely, the connected components of $\hat U^s_{i+1,i}\cap \hat N^s_i$.  Due to the fact that $\partial \hat U^s_{i+1,i}$ contains no heteroclinic point, 
$\partial \hat U^s_{i+1,i}\cap\hat W^{s}_{i,i}$  is isolated from $\hat L^u_i$. Therefore, if the tube 
$\hat N^s_i$ is small enough, $\hat L^u_i$ does not intersect $\partial \hat U^s_{i+1,i}\cap  \hat N^s_i$.
Then, by shrinking $U^s_j,\,j>i+1$ (in the sense of Lemma \ref{Ge}) if necessary, we may guarantee that $\hat U^s_{j,i}\cap \hat N^s_i$  is disjoint from $\partial\hat U^s_{i+1,i}\cap \hat N^s_i$.
 
Since $\hat W^u_{i+2, i}$ accumulates on $\hat W^u_{i+1,i}$, there are only finitely many points of $H_2$ outside of all boxes  $B^1_\ell$, $\ell=1,..., t(1)$.  Let $\bar H_2:= \{h^2_1,...,h^2_{t(2)}\}$ be this finite set. The open set $\hat U^s_{i+2,i}$  is a neighborhood of $\bar H_2$. The  connected  components of $\hat U^s_{i+2,i}\cap \hat N^s_i$ which contain points of $\bar H_2$  will be the box $B^2_\ell$ for $\ell= 1,..., t(2)$. We argue with $B^2_\ell$ with respect to  $\hat L^u_i$ and the neighborhoods $\hat U_{j,i},\, j>i+1,$ in  a similar manner as we do with $B^1_\ell$. And so on, until $\bar H_n$.
 
Due to the induction hypothesis, each above-mentioned box is foliated. Namely, $B^1_\ell$ is foliated by  $\hat F^u_{i+1, i}$; the box $B^2_\ell$ is foliated by $\hat F^u_{i+2, i}$, and so on. But the leaves are not contained in fibres of $\hat N_i$; even more, not every leaf intersects $\hat W^s_{i,i}$.  We have to correct this situation in order to construct the foliation  $F^u_i$ satisfying the requested conditions (i) and (ii).
For every $j>i$, the foliation $F^u_j$ may be extended to the boundary $\partial U^s_j$ and a bit beyond. Once this is done, if $\hat N^s_i$ is  enough shrunk, each leaf of $\hat F^u_{i+k,i}$ {through} $x\in B^k_\ell$
intersects $\hat W^s_{i,i}$ (it is understood that the boxes are intersected with the shrunk tube without changing their names). Thus, we have a projection along the leaves $\pi_{k,\ell}: B^k_\ell\to\hat W^s_{i,i} $; 
but, the image of $\pi_{k,\ell}$ is larger than $B^k_\ell\cap\hat W^s_{i,i}$. Then, we choose a small enlargement $B'^k_\ell$ of $B^k_\ell$ such that   $B'^k_\ell\smallsetminus B^k_\ell$ is foliated by $\hat F^u_{i+k,i}$ and avoids  the lamination ${\hat L^u_i}$. On $B'^k_\ell\smallsetminus B^k_\ell$ we have two projections: one is $\hat\pi^u_i$ and the other one is $\pi_{k,\ell}$. We are going to interpolate between both using a partition of unity (we do it for $B^k_\ell$ but it is understood that it is done for all boxes). Let $\phi:\hat N^{s}_{i}\to[0,1]$ be a smooth function which equals 1 near   $B^k_\ell$ and whose support is contained in $B'^k_\ell$. Define a global $C^1$ retraction $\hat q: \hat N^{s}_{i}\to \hat W^s_{i,i}$ by the formula $$\hat q(x)=\bigl(1-\phi(x)\bigr)  \hat\pi^u_i(x) + \phi(x)\bigl(\pi_{k,\ell}(x)\bigr).$$
Here, we use an affine manifold structure on each component of $\hat W^s_{i,i}$ by identifying it with the 1-torus $\mathbb T:= \mathbb R/\mathbb Z$. So, any positively weighted barycentric combination makes sense for a pair of points sufficiently close. When $x\in \hat W^s_{i,i}$, we have $\hat q(x)=x$. Then, by shrinking  the tube $\hat N^{s}_{i}$ once more if necessary we make  $\hat q$ be a fibration whose fibres are transverse to the lamination $\hat L^s_i$ and we make each leaf of $\hat F^u_{j,i}, \,j>i,$ in every  box $B^\ell_k$ be contained in a fibre of $q$. Henceforth, 
taking the preimage of that  tube (and its fibration) by $p_i$  and adding the unstable manifold $W^u_i$ provide the requested $U^s_i$ and its foliation $F^u_i$ satisfying the required properties. Thus, the induction is proved. 
\end{demo}

We also have the following statement.
\begin{lemm} \label{F^s_i}
There exist  $f$-invariant neighborhoods $U^{u}_{0},\dots,U^{u}_{n}$ of  the sets $\Si_{0},\dots,\Si_{n}$  respectively, equipped with one-dimensional $f$-invariant foliations ${F}^{s}_{0},\dots,{F}^{s}_{n}$ {of class $C^{1,0}$} such that the  following properties  hold for each $i\in\{0,\dots,n\}$:

{\rm (iii)} the stable manifold $W^s_i$ is a leaf of the foliation ${F}^{s}_{i}$ and each leaf of the foliation ${F}^{s}_{i}$ is transverse to $L^u_i$;
 
{\rm (iv)} for any $0\leq j<i$ and $x\in U^u_{i}\cap U^u_{j}$, we have the inclusion $({F}^s_{j,x} \cap{U}^u_{i})\subset {F}^s_{i,x}$.
\end{lemm}
\begin{demo} The proof is done by an increasing induction from $i=0$; it is skipped due to similarity to the previous one.
\end{demo}
\medskip

Before entering Step 2,  we recall the definition of fundamental domain for a free action.
\begin{defi} Let $g:X\to X$ be a homeomorphism acting freely on $X$. A closed subset $D\subset X$ is said to be a
\emph{fundamental domain} for the action of $g$ if the following properties hold:

{\rm 1.} $D$ is the closure of its interior $\mathop{D}\limits^\circ$;

{\rm 2.} $g^k(\mathop{D}\limits^\circ)\cap D=\emptyset$ for every integer $k\neq 0$;

{\rm 3.} $X$ is the union $\cup_{k\in \Z} \, g^k(D)$.

\end{defi} 

{\bf Step 2.} We prove the following statement for each $i=0,\dots,n$. 
\begin{lemm}\label{fund} $ $
{\rm (v)} There exists  an $f$-invariant neighborhood $\tilde N_{i}$ of the set $\Si_i$ contained in $U^s_{i}\cap U^u_{i}$ and such that the restrictions of the foliations ${F}^u_i$ and ${F}^s_i$ to $\tilde N_i$ are transverse.
\end{lemm}

\begin{demo} For this aim, let us choose a fundamental domain
$K^s_i$ of the restriction of $f$ to $W^s_i\smallsetminus\Si_i$ and take a tubular neighborhood $N(K^s_i)$ of  $K^s_i$ whose disc fibres are contained in leaves of $F^u_i$. {By construction, $F^u_i$ is transverse to $W^s_i$ and, according to the Lemma \ref{F^s_i}, $F^s_i$ is a $C^{1,0}$-foliation. Therefore, if the tube $N(K^s_i)$ is small enough, $F^u_i$ is transverse to $F^s_i$ in $N(K^s_i)$.}  
Set $$\tilde N_i:=W^u_i\,\bigcup_{k\in\mathbb Z}f^k\!\left(N(K^s_i)\right). $$ 
This is a neighborhood  of $\Si_i$\,; it satisfies condition (v) and the previous properties (i)--(iv) still hold.  A priori the boundary of $\tilde N_i$ is only piecewise smooth; but, by choosing  $N(K^s_i)$ correctly at its corners we may arrange that $\partial\tilde N_i$ be smooth.
\end{demo}

{\bf Step 3.} {For proving Theorem \ref{compatible}} it remains to show the existence of linearizable 
neighborhoods $N_i\subset\tilde N_i,~i=0,\dots,n$, for which the required foliations are the restriction to 
$N_i$ of the foliations {${F}^u_i$ and ${F}^s_i$.} For each orbit of $f$ in $\Si_i$, choose one $p$.
 Let $\tilde N_p$ be a connected component of $\tilde N_i$ containing $p$. There is a homeomorphism 
 $\varphi^u_p:W^u_p\to W^u_{O}$ (resp. $\varphi^s_p:W^s_p\to W^s_{O}$) conjugating the 
 diffeomorphisms $f^{per(p)}|_{W^u_p}$ and $\mathcal E|_{W^u_{O}}$ (resp. $f^{per(p)}|_{W^s_p}$ and 
 $\mathcal E|_{W^s_{O}}$). In addition, for any point $z\in\tilde N_p$ there is unique pair of points
  $z_s\in W^s_p,~z_u\in W^u_p$ such that  $z=F^s_{i,z_u}\cap F^u_{i,z_s}$. We define a topological 
  embedding $\tilde{\mu}_{p}:\tilde N_{p}\to\mathbb R^3$ by the formula $\tilde{\mu}_{p}(z)=(x_1,x_2,x_3)$ 
  where  
$(x_1,x_2)={\varphi}^u_{p}(z_u)$ and $x_3={\varphi}^s_{p}(z_s)$. 
Since the foliations $F^u_i$ and $F^s_i$ are $f$-invariant, this definition makes $\tilde\mu_p$ conjugate 
the restriction $f^{per(p)}\vert_{\tilde N^p}$ to $\mathcal E^{per(p)}$. For $k= 1, \dots, per(p)-1$, set 
$ \tilde N_{f^k(p)}:=f^k(\tilde N_p)$ and define $\tilde\mu_{f^k(p)}$ so that the equivariance formula holds:
$\tilde\mu_{f^k(p)}(f^k(x))=a^k \tilde\mu_{p}(x)$ for every $x\in \tilde N_{p}$. Choose $t_0\in(0,1]$ such
 that $\mathcal N^{t_0}\subset\tilde\mu_p(\tilde{N}_{p})$ {for every $p\in\Si_i$.} Observe that 
 $\mathcal E\vert_{\mathcal{N}^{t_0}}$ is conjugate to 
 $\mathcal E\vert_{\mathcal{N}}$ by the suitable homothety $h$.  Set $N_{p}=\tilde\mu^{-1}_p(\mathcal{N}^{t_0})$ 
 and $\mu_{p}= h \tilde\mu_p:N_{p}\to \mathcal{N}$. Then, $N_p$ is the requested neighborhood with its 
 linearizing homeomorphism $\mu_{p}$. This finishes the proof of Theorem \ref{compatible}.
\end{demo}

\section{Proof of the classification theorem}\label{IV}
Let us prove that {the diffeomorphisms $f$ and $f'$ in}  $G(M)$ are topologically conjugate if and only if there is a homeomorphism $\hat\varphi:\hat V_f\to\hat  V_{f'}$ such that 

$(1)\,\eta_{_f}=\eta_{_{f'}}\hat\varphi_*$;

$(2)\,\hat\varphi(\hat{{\Ga}}^s_{f})=\hat{{\Ga}}^s_{f'}$ and  $\hat\varphi(\hat{\Ga}^u_{f})=\hat{\Ga}^u_{f'}$.
\subsection{Necessity} 
{Let $f:M\to M$ and $f':M\to M$ be two elements in $G(M)$ which are  topologically conjugated by some} homeomorphism $h:M\to M$. Then $h$ conjugates the invariant manifolds of periodic points of $f$ and $f'$. More precisely, if $p$ is a periodic point of $f$, then $h(p)$ is a periodic point of
$f'$ and $h\left(W^{u}(p)\right)= W^{u}(h(p))$, $h\left(W^{s}(p)\right)= W^{s}(h(p))$. In particular, $h$ maps $V_f$ to $V_{f'}$ by a homeomorphism  noted $\vp$. Moreover, if $x$ is any points of $V_f$, for every $n\in \mathbb Z$ the following holds:
$$\vp(f^n(x))=f'^n(\vp(x)).$$
This formula says exactly that $\vp$ is the lift of a map $\hat\vp: \hat V_f\to \hat V_{f'}$. By construction of $\eta_f$, the same formula says that $\eta_f= \eta_{f'}\circ \hat\vp_*$, where $ \hat\vp_*:H_1(\hat V_f; \mathbb Z)\to H_1(\hat V_{f'}; \mathbb Z)$ denotes the map induced in homology. By definition of the quotient topology, $\hat\vp$ is continuous. Since the same holds for $\vp^{-1}$, one checks that $\hat\vp$ is a homeomorphism. 
As $\vp$ conjugates the laminations {$\Ga^s_f$} (resp. $\Ga^u_f$) to 
$\Ga^{s}_{f'}$ (resp. $\Ga^{u}_{f'}$), the same holds for  $\hat\vp$ in the quotient spaces with respect the projections of the laminations.

\subsection{Sufficiency} 
For proving the sufficiency of the conditions in Theorem \ref{t.invariant}, let us consider  a homeomorphism $\hat\varphi:\hat V_{f}\to\hat V_{f'}$ such that:

(1) $\eta_{_f}=\eta_{_{f'}}\hat\varphi_*$;

(2)  $\hat\varphi(\hat{\Gamma}^s_{f})=\hat{\Gamma}^s_{f'}$ and  $\hat\varphi(\hat{\Gamma}^u_{f})=\hat{\Gamma}^u_{f'}$.

\nd From now on, the dynamical objects attached to $f'$ will be denoted by {$L'^u, L'^s, \Si'_i, \ldots$ with 
the same meaning as  $L^u, L^s, \Si_i, \ldots$} have with respect to  $f$. By (1), $\hat\vp$ lifts
 to an {\it equivariant} homeomorphism $\varphi: V_{f}\to V_{f'}$,  that is: $f'|_{V_{f'}}= \vp f \vp^{-1}|_{V_{f'}}$  (for brevity, equivariance stands for $(f,f')$-equivariance). By (2), $\varphi$ maps  $ \Gamma_f^u$ to $\Gamma_{f'}^u$ and  $ \Gamma_f^s$ to $\Gamma_{f'}^s$.
Thanks to Theorem \ref{compatible} we may use compatible linearizable neighborhoods of the saddle points of $f$ (resp. $f'$). 

An idea of the proof is the following: we  modify the homeomorphism $\varphi$ in a neighborhood of $\Ga_f^u$ such that {the} final homeomorphism preserves {the} compatible foliations, then we 
do  similar modification near $\Ga^s_f$. So we get a homeomorphism $h:M\setminus(\Om_0\cup\Om_3)\to M\setminus(\Om'_0\cup\Om'_3)$ conjugating $f|_{M\setminus(\Om_0\cup\Om_3)}$ with $f'|_{M\setminus(\Om'_0\cup\Om'_3)}$. Notice that $M\setminus(W^s_{\Om_1}\cup W^s_{\Om_2}\cup\Om_3)=W^s_{\Om_0}$ and $M\setminus(W^s_{\Om'_1}\cup W^s_{\Om'_2}\cup\Om'_3)=W^s_{\Om'_0}$. Since $h(W^s_{\Om_1})=W^s_{\Om'_1}$ and 
$h(W^s_{\Om_2})=W^s_{\Om'_2}$ then $h(W^s_{\Om_0}\setminus\Om_0)=W^s_{\Om'_0}\setminus\Om'_0$. Thus for each connected component $A$ of $W^s_{\Om_0}\setminus\Om_0$ there is a sink $\omega\in\Om_0$ such that $A=W^s_\omega\setminus\omega$. Similarly $h(A)$ is a connected component of $W^s_{\Om'_0}\setminus\Om'_0$ such that  $h(A)=W^s_{\omega'}\setminus\omega'$ for a sink $\omega'\in\Om'_0$. Then we can continuously extend $h$ to $\Om_0$ assuming $h(\omega)=\omega'$ for every $\omega\in\Om_0$. A similar extension of $h$ to $\Om_3$ finishes the proof. Thus below in a sequence of lemmas we explain only how to modify the homeomorphism $\varphi$ in a neighborhood of $\Ga_f^u$ such that {the} final homeomorphism preserves {the} compatible foliations.

Recall the partition $\Si_0\sqcup\cdots\sqcup \Si_n$ associated with the  Smale order 
on the periodic points {of index 2}.

\begin{lemm}\label{leaf_number} For every $i= 0, \dots, n$ the following equality holds 
 $\vp(W^u_i\cap{V_f})=W^{\prime u}_i\cap{V_{f'}}$ and there is a unique continuous extension 
{of} $\vp|_{W^u_i\cap V_f}$ to $\Si_i$ which is equivariant and bijective from $\Si_i$ to $\Si'_i$.
\end{lemm}
\begin{demo} {Let $p\in \Si_0$. Denote its orbit by $orb_f(p)$. The punctured unstable manifold 
$W^u(p)\smallsetminus\{p\}$ projects by $p_f$ to one compact leaf $\ell(p)$.  Both sides 
of the next equality are $f$-invariant and project to the same leaf, thus:
$$p_f^{-1}(\ell(p))= W^u(orb(p))\smallsetminus\{orb(p)\}.$$
Then, the number of connected components of $p_f^{-1}(\ell(p))$ is $per(p)$, the period of $p$. The 
image $\hat\vp(\ell(p))$ is a compact leaf of $\hat\Ga^u_{f'}$. By the previous argument,  it is $\ell'(p')$ for 
some $p'\in \Si'_0$. Since {$\hat\vp$ lifts to $\vp$,} then 
$\vp\left[ p_f^{-1}(\ell(p))\right]=  p_{f'}^{-1}(\ell'(p'))$ 
which implies the equality of the number of connected components. Thus $per(p)=per(p')$. 
From this, we can deduce that, up to replacing $p'$ with $f'^k(p')$ for some integer $k$,
we have $\vp\left(W^u(p)\smallsetminus \{p\}\right)= W^u(p')\smallsetminus \{p'\}$. Using the property 
$p= \lim\limits_{n\to -\infty} f^n(x)$ for every $x\in W^u(p)$ and the similar property for $p'$ in addition to  
the  equivariance of $\vp$, one extends continuously $\vp|_{W^u_p}$ by defining $\vp(p)=p'$. Doing 
the same for every orbit of $ \Si_0$, we get a continuous extension of {$\vp|_{W^u_0}$}  to $\Si_0$ 
which is still equivariant. One easily {checks} that this extension is {continuous, unique, and hence equivariant}. Then, arguing similarly with $\hat\vp^{-1}$, we derive that the extension of $\vp$ maps $\Si_0$ bijectively onto $\Si'_0$. 
 
Denote $\ell_0:= \bigcup\limits_{p\in \Si_0}\ell(p)$. We have $\hat\vp(\ell_0)=\ell'_0$. Let now $p\in \Si_1$.
 The closure in $\hat V_f$ of $\ell(p):=p_f(W^u(p)\smallsetminus \{p\})$ is contained in $\ell(p)\cup\ell_0$. 
 We deduce that $\hat\vp(\ell(p))$ is a leaf of $\hat\Ga^u_{f'}$ of the form $\ell'(p')$ for some $p'\in \Si'_1$ 
 and we can continue inductively. Thus, there is a continuous extension of $\vp|_{W^u_i}$ to every
  $\Si_i$ for $i= 0, 1,\dots, n$ which is a bijection $\Si_i\to\Si'_i$. Arguing with $\hat\vp^{-1}$, we derive 
  that $n'=n$.}
\end{demo}

Recall the radial functions $r^u_i, r^s_i: N_i\to [0,+\infty)$ which are introduced above Definition \ref{2dopsystem}; recall also the order which is defined by $r^s_i$ on each stable separatrix  {$\ga_p$} of $p\in \Si_i$. Analogous functions are associated with the dynamics of $f'$.
  
\begin{lemm} \label{pusk} There is a unique continuous extension of $\vp|_{\Ga^u_f}$
$$\vp^{us}:\Ga^u_f\cup \left(L^u\cap L^s\right) \longrightarrow \Ga^u_{f'}\cup \left(L'^u\cap L'^s\right)$$
such that {the following holds}:
 
\nd {\rm (1)} If  $x\in W^u_j\cap W^s_i$, $j>i$, then $\vp^{us}(x)\in  W'^u_j\cap W'^s_i$.
  
\nd {\rm (2)} If $x$ and $y$ lie in $\gamma_p\cap L^u$ with $r^s_p(x)<r^s_p(y)$,  then $\vp^{us}(x)$ and 
$\vp^{us}(y)$ lie in $ \gamma'_{\varphi(p)}\cap L'^u$ with
 ${r'^s_{\varphi(p)}}(\vp^{us}(x))< {r'^s_{\varphi(p)}}(\vp^{us}(y))$.
\end{lemm}

Notice that $\vp|_{\Ga^u_f}$ being equivariant, its continuous extension is also equivariant.\\

\nd\begin{demo} This statement is proved by induction on $i$. We recall that 
$V_i\smallsetminus L^s_i=V_f\smallsetminus cl(W^u_i)$ is a dense  open set in $V_f$ (and similarly  with 
$^\prime$), and according to Lemma \ref{leaf_number}, $\vp$ maps $V_i\smallsetminus L^s_i$ to 
$V'_i\smallsetminus L'^s_i$ homeomorphically  and conjugates $f$ to $f'$. Thus, for every $i=0,\dots, n$,
we have an equivariant homeomorphism $\vp_i: V_i\smallsetminus  L^s_{i}\to  V'_i\smallsetminus L'^s_{i}$ 
which maps $W^u_{j}\smallsetminus L^s_{i}$ to $W'^u_{j}\smallsetminus L'^s_{i}$ for every $j>i$, again
 as a consequence of Lemma \ref{leaf_number}.
 
First, take $i=0$.  The manifold $\hat V_0$ is closed and three-dimensional. We have $\hat L^s_0= \hat W^s_{0,0}$\,, which consists of finite number disjoint  smooth circles, and similarly with $^\prime$. 
\begin{figure}[h]
\centerline{\includegraphics[width=13cm,height=9cm]{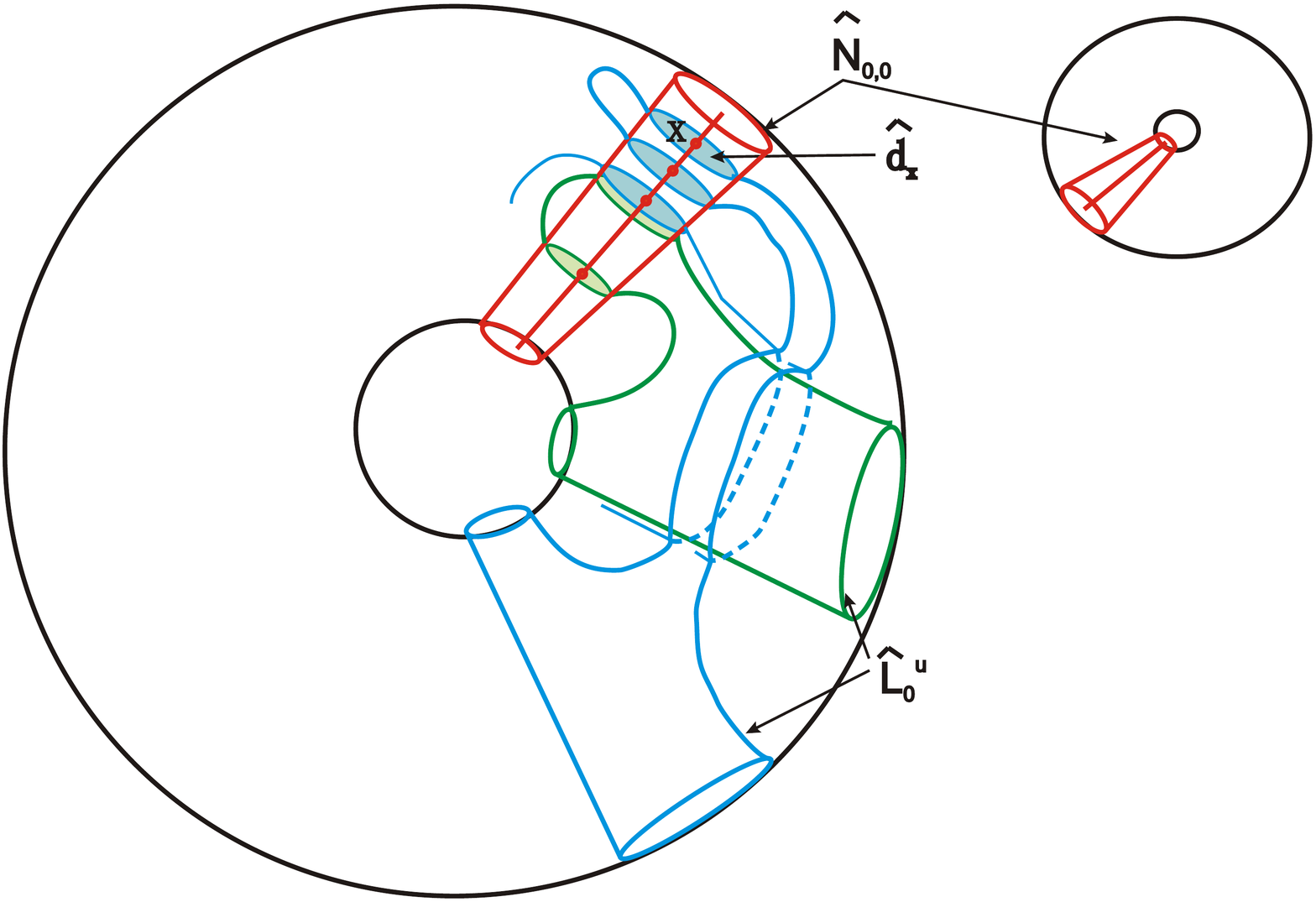}}\caption{{Case $i=0$ in proof of Lemma \ref{pusk} for the} diffeomorphism from Figure \ref{faz}.} \label{dx}
\end{figure}
   
We look for an extension $\hat\vp^{us}_0$ of $\hat\vp_0|_{\hat L^u_0}$ to $\hat L^s_0\cap \hat L^u_0$. 
If $N_0$ is the neighborhood {of $\Si_0$ extracted  from a compatible system given by Theorem 
\ref{compatible} and if} $\hat N_{0,0}$ denotes  the corresponding tubular neighborhood of $\hat L^s_0$ 
in $\hat V_0$, the trace of $\hat L^u_0$ in that tube is a lamination by disks:
$$\hat L^u_0\cap \hat N_{0,0}=\{ \hat d_x\mid x\in \hat L^u_0\cap \hat L^s_0\},$$ where $\hat d_x$ denotes 
the fiber of the tube over $x\in \hat L^s_0$ (see figure \ref{dx}). 

{The complement in $\hat V'_0$ of the interior of $\hat N'_{0,0}$ is a compact set contained in
$\hat V'_0\smallsetminus \hat L'^s_0$. Then its preimage $K$ by the homeomorphism $\hat\vp_0$
is a compact set contained in $\hat V_0\smallsetminus \hat L^s_0$. When $t=0$, we have 
$\hat N^t_{0,0}= \hat L^s_0$ and hence disjoint from $K$. Then, if $t$ is small enough,
$\hat\vp_0(\hat N^t_{0,0}\smallsetminus \hat L^s_0)\,\subset\,\hat N'_{0,0}\smallsetminus \hat L'^s_0 $.
 Finally,} the map
 $\hat \vp_0$ (which is not defined on $\hat L^s_0$) possesses the two following properties:
\begin{enumerate}
\item If $N_0$ is shrunk enough, we have $(\hat\vp_0(\hat N_{0,0})\cap\hat L^u_0)\subset (\hat N'_{0,0}\cap\hat L^{\prime u}_0)$\,, where $\hat N'_{0,0}$ denotes the tube associated with the chosen linearizable neighborhood $N'_0$ of {$\Si'_0$}.
\item If $\hat d_x$ is a plaque of $\hat L^u_0\cap \hat N_{0,0}$\,, the image $\hat\vp_0(\hat d_x\smallsetminus \{x\})$ is contained in some fiber $\hat d_{x'}$, with $x'\in \hat L'^u_0\cap L'^s_0$. 
\end{enumerate}
As a consequence, the requested extension may be defined by $\hat\vp^{us}_0(x)=x'$.  As the considered 
plaques are arcwise connected, the construction lifts to the cover and yields a continuous map 
$\vp^{us}_0: \Ga^u_f\cup (L^s_0\cap L^u_0)\to \Ga^u_{f'}\cup(L'^s_0\cap L'^u_0)$ 
which is a continuous {equivariant} extension of $\vp|_{\Ga^u_f}$. 
  
It remains to prove that $\vp^{us}_0$ is increasing on its domain in each separatrix of {$\Si_0$}. For this aim, consider a point {$p\in\Si_0$, one of its separatrices $\ga_p$}  and a connected component $N_{\gamma_p}$ of $N_p\smallsetminus W^u_p$ containing $\gamma_p$. Take an infinite proper arc 
$C$ in $N_{\gamma_p}\smallsetminus W^s_p$ which crosses  transversely each leaf of the foliation $F^u_0$ and 
 which has one end in $p$. 
 We orient $C$ so that its projection onto $\ga_p$ is positive.
Its image through $\vp_0$ is a proper arc $C'$ contained in  {$N'_{\varphi(p)}\smallsetminus W'^u_{\varphi(p)}$.}  Moreover, ${\varphi(p)}$ is  one  end of $C'$. 
For  $x,y\in \gamma_p\cap L^u$, the inequality $r^s_p(x)<r^s_p(y)$ implies
${r'^s_{\varphi(p)}}(\vp^{us}_0(x))< {r'^s_{\varphi(p)}}(\vp^{us}_0(y))$ if we are sure that $C'$ intersects each 
leaf of $L'^u_0\cap N'_{\varphi(p)}$ at most in one point. 
That is true since $\vp_0$ is a homeomorphism on {its} image from $N_0\smallsetminus W^s_0$ to $N'_0\smallsetminus W'^s_0$ mapping $L^u_0$ into $L'^u_0$. 
  
\begin{figure}[h]
\centerline{\includegraphics[width=13cm,height=9cm]{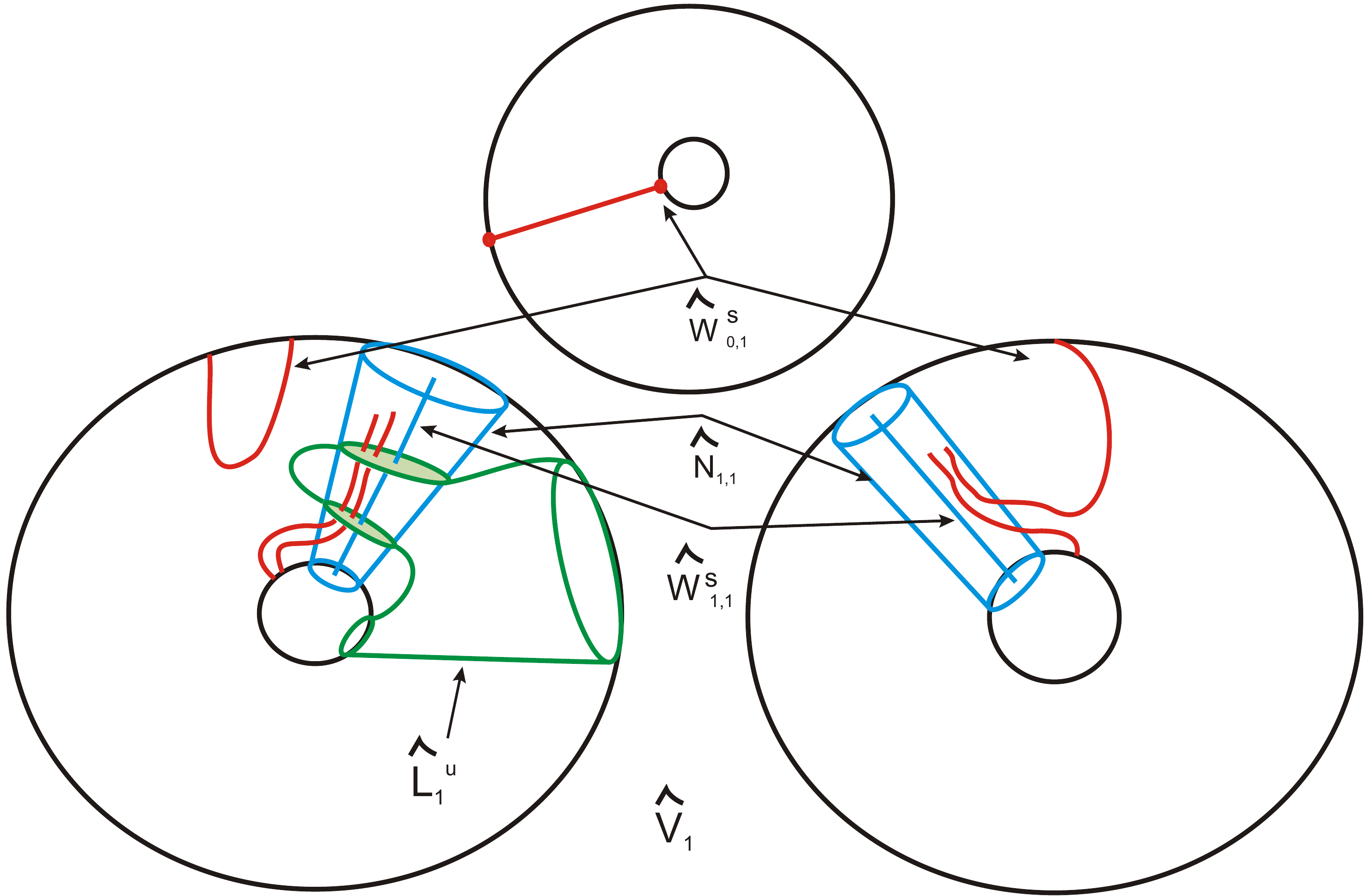}}\caption{Illustration {of induction in Lemma \ref{pusk} for the diffeomorphism} from Figure \ref{faz}.} \label{induc}
\end{figure}
For the induction,  {let $i\in\{1,\ldots,n\}$ and let us assume that there is  a continuous} extension 
$${\vp^{us}_{i-1}}: \Ga^u_f \cup\bigcup\limits_{j=0}^{i-1}(L^s_j\cap L^u_j)\to \Ga^u_{f'} \cup\bigcup\limits_{j=0}^{i-1}(L'^s_j\cap L'^u_j),$$
{which is monotone on each separatrix of $\Si_j,\ j<i$.}
{The image $\hat W^s_{i,i}$ of $W^s_i$ by the projection $p_i:V_i\to\hat V_i$
is made of finitely many disjoint circles which are the images of the stable 
 separatrices of $\Si_i$. About $\hat W^s_{i,i}$\,, there is a tube $\hat N_{i,i}$ which is the projection 
 by $p_i$ of a neighborhood $N_i$ of $\Si_i$ extracted  from a compatible system given by Theorem \ref{compatible} (see Figure \ref{induc}).
The trace of $\hat L^u_i$ in that tube is a lamination by disks:
  $$\hat L^u_i\cap \hat N_{i,i}=\{\hat d_x\mid x\in \hat L^u_i\cap \hat W^s_{i,i}\},$$ where $\hat d_x$ denotes the fiber of the tube over $x\in \hat W^s_{i,i}$.
  
 In $V_i$, there are two laminations $L^u_i$ and $L^s_i$ (and the corresponding objects with $'$).
 The map $\vp_i$, not defined on $L^s_i$, sends $L^u_i\smallsetminus L^s_i$
 homeomorphically onto $L'^u_i\smallsetminus L'^s_i$. By the induction hypothesis, 
 the restriction $\vp_i \vert (L^u_i\smallsetminus L^s_i)$ extends continuously 
 to $L^u_i\smallsetminus W^s_i$; this extension, automatically equivariant, is denoted by $\psi_i$.
 This induces on the quotient space $\hat V_i$ a homeomorphism 
 $$\hat\psi_i: \hat L^u_i\smallsetminus \hat W^s_{i,i} \to \hat L'^u_i\smallsetminus \hat W'^s_{i,i}\,.$$
 
 In order to extend $\hat\psi_i$ to $\hat L^u_i\cap \hat L^s_i$, we use the fact that $\hat L^u_i$
 is compact   for arguing as in the case $i=0$. Consider the tube $\hat N^t_{i,i}$ depending on $t\in(0,1)$
  and look at its compact lamination by disks $\hat L^u_i\cap \hat N^t_{i,i}$. After removing $\hat W^s_{i,i}$
  which marks one puncture on each leaf,
   it is leaf-wise mapped by $\hat \psi_i$ into $\hat V'_i\smallsetminus \hat W'^s_{i,i} $. As in case $i=0$,
the above-mentioned compactness allows us to conclude that there exists some $t\in(0,1)$ such that
   $\hat L^u_i\cap( \hat N^t_{i,i}\smallsetminus \hat W^s_{i,i})$ is mapped into $\hat N'_{i,i}$
   where $\hat N'_{i,i}$ denotes the tube associated with the chosen linearizable neighborhood $N'_i$
   of $\Si'_i$. Finally,}
 the map $\hat \psi_i$ possesses {the} two following properties: 
 
\begin{enumerate}
\item[1.] If $N_i$ is shrunk enough, we have $(\hat\psi_i(\hat N_{i,i})\cap\hat L^u_i)\subset (\hat N'_{i,i}\cap\hat L^{\prime u}_i)$\,.
\item[2.] If $\hat d_x$ is a plaque of $\hat L^u_i\cap \hat N_{i,i}$\,, the image $\hat\psi_i(\hat d_x\smallsetminus \{x\})$ is contained in some fiber $\hat d_{x'}$, with $x'\in \hat L'^u_i\cap W'^s_i$. 
\end{enumerate}
Now, the extension $\hat\vp^{us}_i$ of $\hat \psi_i$ is defined by $x\mapsto x'$. One checks it is a 
continuous extension. The requested $\vp^{us}_i$ is the lift of $\hat\vp^{us}_i$ to $V_i$. It has the required 
properties allowing us to finish the induction.
\end{demo}
\begin{remark} \label{mai} {\rm Due to Lemma \ref{Ge} we may assume that in all lemmas below 
{the chosen values  $t=\beta_i,a_i,...$  are} such that the boundary of the linearizable neighborhood 
 $N^t_i$ does not contain any heteroclinic point.} 
\end{remark} 
  
\begin{lemm} \label{psibe} There are numbers $\beta_0,\dots,\beta_n\in (0,1)$  
such that, for  every {$i\in\{0,\dots,n\}$, for every point $p\in\Si_i$ and $x\in N^{\beta_i}_p\cap L^u$, 
 the following inequality holds:
  $$r'^u_{i}(\vp^{us}(x^u_i)) r'^s_{i}(\vp^{us}(x^s_i))<\frac12\,.$$} 
\end{lemm}
\nd\begin{demo} As $N_n\cap L^u=W^u_n$ and $\vp^{us}(W^u_n)=W'^u_n$, it is possible to chose 
any $\beta_n\in(0,1)$.
 
Indeed, for  $p\in \Si_n$ and $x\in W^u_p$, we have $r'^s_{\vp(p)}(\vp^{us}
(x^s_i))=0$.
  For $i\in\{0,\dots,n-1\}$ {and $p\in\Si_i$,  choose some heteroclinic point $y\in W^s_p\cap L^u$ arbitrarily. 
  Set:
$$\la'^u_p(t)=\sup_x\{r'^u_{\vp(p)}(\vp^{us}(x^u_i))\mid  x\in N^t_p\cap F^u_{i,y}\} 
\quad{\rm and}\quad
 \la'^s_p=r'^s_p(\vp^{us}(y))\,.$$}
{When $t$ goes to $0$, the arc $N^t_p\cap F^u_{i,y}$ shrinks to the point $y$. Then, according to Lemma \ref{pusk}, $\la'^u_p(t)$ also goes to $0$. Therefore,} there exists some $\beta_p\in(0,1)$ such that
  $\la'^u_p(\beta_p)\la'^s_p<\frac18$. Denote by $Q_p$ {the compact subset of $M$ bounded by 
  $\partial N^{\beta_p}_p, F^u_{i,y}$ and $f^{per(p)}(F^u_{i,y})$.} Notice that $Q_p$ is a 
  fundamental domain {for the restriction of $f^{per(p)}$ to the} connected component of 
  $N^{\beta_p}_p\smallsetminus W^u_p$ containing $y$. For {every} $x\in Q_p$, 
we have $r'^u_{\vp(p)}(\vp^{us}(x^u_i))\leq 4\la'^u_p(\beta_p)$ and 
$r'^s_{\vp(p)}(\vp^{us}(x^s_i))\leq \la'^s_p$. 
 Then, for every  $x\in Q_p\cap L^u$ we have:
$$r'^u_p(\vp^{us}(x^u_i))\, r'^s_p(\vp^{us}(x^s_i))\leq 
4\la'^u_p(\beta_p)\,\la'^s_p<\frac12.$$
Set $\beta_i=\min\limits_{p\in\Si_i}\{\beta_p\}$. {Hence, $\beta_i$ is the required number.} 
\end{demo}

\begin{lemm} \label{psu2} When $n>0$, there exist real numbers $a_j\in(0,\beta_j]$ fulfilling 
the following property: for every $j= 1,\ldots, n$ and every integer $i<j$, each connected component
of $\hat W^s_{i,i} \cap \hat N^{a_j}_{j,i}$ is an open interval which is either disjoint from 
$A_j^i:=\bigcup\limits_{k=i+1}^{j-1}\hat N^{a_k}_{k,i}$ or included in $A_j^i$. Moreover,
only finitely many of these intervals are not covered by $A_j^i$.
\end{lemm}
\begin{demo} The proof is done by induction on $j$ from 1 to $n$. For  $j=1$, one is allowed to take 
$a_1= \beta_1$. Indeed, $\hat W^s_{0,0}$ is a  smooth curve and $\hat W^u_{1,0}$ is a smooth closed
surface which is transverse to $\hat W^s_{0,0}$. Therefore, there are finitely many intersection points.
 By the choice of $\beta_1$, the projection in $\hat V_0$
of $N_1^{a_1}$ is a tubular neighborhood of $\hat W^u_{1,0}$. Moreover, each component of 
$\hat W^s_{0,0}\cap \hat N_{1,0}^{a_1} $ is a fiber of this tube.

For the induction, assume the numbers $a_1, \ldots, a_{j-1}$ are given with the required properties and let 
us find $a_j$. In particular, the subset $A_j^i$ is assumed to be defined. According to Remark \ref{mai},
the boundary of $A_j^i$ contains no heteroclinic point.

First, fix $i<j$. Consider the projection $\hat W^u_{j, i}$ of $W^u_j$ in $\hat V_i$. This is a union of leaves
 in the lamination $\hat L^u_i$. The following is a well-known fact (see, for example, Statement 1.1 in \cite{GrPo2013}): if $x$ is a point from $\hat L^u_i$ which is accumulated
by a sequence of plaques  from $\hat W^u_{j,i}$, then $x$ does not lie in $\hat W^u_{j,i}$ but 
belongs to some $\hat W^u_{k,i}$ with $k<j$. Then the part of $\hat W^u_{j,i}$ which is covered 
by $A_j^i$  contains every intersection points $\hat W^s_{i,i}\cap \hat W^u_{j,i}$
except  finitely many of them. From this finiteness and the fact that 
$A_j^i\cap\hat W^s_{i,i}\cap \hat W^u_{j,i}$ is actually contained in $A_j^i$, an easy 
compactness argument allows us to find a positive number $a_j^i$ such that the collection of 
disjoint intervals made by $\hat W^s_{i,i}\cap \hat N_j^{a_j^i}$ fulfills the requested property
with respect to the considered $i$. Indeed, let  $B(t)$ be the closure of $\partial A_j^i\cap\hat W^s_{i,i}\cap \hat N_j^{t} $. The intersection $\bigcap\limits_{k\in \mathbb N}B(\frac1k)$ is empty. Then $B(t)$ is empty when $t$ is small enough.
 
By defining $a_j:= \inf\{a_j^0,\ldots, a_j^{j-1}\}$, we are sure 
that $\hat N_j^{a_j}$ satisfies all the requested properties.
\end{demo}

The corollary below  follows  from {Lemma} \ref{psu2} immediately. 
\begin{collary} \label{con} For each $i\in\{0,\dots,n-1\}$ the intersection
$\hat W^{s}_{i,i}\cap (\bigcup\limits_{j=i+1}^n\hat{N}^{a_{j}}_{j,i})$ consists of 
finitely many open arcs $\hat I^{i}_1,\dots,\hat I^{i}_{r_{i}}$ such that{, for each $l=1,\dots,r_{i},$ the arc $\hat I^{i}_l$}  is a connected component of 
 $\hat W^{s}_{i,i}\cap\hat{N}^{a_{j}}_{j,i}$ for some $j>i$.
\end{collary}

For brevity, for $i=0,\ldots, n$, we denote by $\vp^{u}_i$ the restriction $\vp^{us}|_{W^u_i}$ 
in the rest of the proof of Theorem \ref{t.invariant}. Let {$\psi_i^s: W^s_i\to W'^s_i$} be any equivariant homeomorphism which extends {$\vp^{us}\vert_{W_i^s\cap L^u}$} and let   $t_i\in (0,1)$ be a small enough number so that, for every $x\in N_i^{t_i}$, the next inequality holds:
$$(*)_i\quad\quad r'^s(\vp_i^u(x^u_i)) \,r'^u(\psi^s_i(x^s_i))<1.$$
In this setting, one derives an equivariant embedding $\phi_{\vp^u_i ,\psi^s_i}: N_i^{t_i} \to N'_i$ which is defined by sending $x\in N_i^{t_i}$ to $\left(\vp_i^u(x^u_i),\psi^s_i(x^s_i)\right)$.

\begin{lemm} \label{1dim} There is an equivariant homeomorphism $\psi^s:L^s\to L'^{s}$ consisting of conjugating  homeomorphisms $\psi^s_{0}:W^s_{0}\to W^{\prime s}_{0},\dots,\psi^s_{n}:W^s_{n}\to W^{\prime s}_{n}$ such that for each $i\in\{0,\dots,n\}$:
\begin{itemize}
\item[\rm (1)] $\psi^s_{i}|_{W^s_{i}\cap L^u}={\vp^{u}_i}|_{W^s_i\cap L^u}$;
\item[\rm (2)] the topological embedding ${\phi}_{{\vp}^{u}_i,\psi^s_i}$ is well-defined on $N^{a_i}_i$;
\item[\rm (3)] if $x\in (W^s_i\cap N_j^{a_j}),\ j>i$, then $\psi^s_i(x)={\phi}_{{\vp}^{u}_j,\psi^s_j}(x) $.
\end{itemize}
\end{lemm}
\nd\begin{demo} We are going to construct $\psi^s_i$ by a decreasing induction on $i$ from $i=n$  {to}
 $i=0$. The stable manifolds of the saddles in $\Si_n$  have no heteroclinic points. Therefore,  the only 
constraints on $\psi^s_n$ imposed by the first item is its value on  $\Si_n$. In particular, we are allowed to 
change $\psi^s_n$ to $f'^k\circ \psi^s_n$ if $k$ is {\it admissible} in the sense that $k$ is 
a multiple of all periods $per(p), \ p\in \Si_n$. 

This remark is used in the following way. One starts with any equivariant homeomorphism 
$\psi^s_n$ such that for any $p\in \Si_n$ the stable manifold $W^s_p$ is mapped to the stable manifold 
of $\vp_n^u(p)$; hence, item 1 is fulfilled. Choose a fundamental domain $I$ of 
$f\vert_{W^s_n\smallsetminus\Si_n}$. Consider the fundamental domain of 
$f\vert_{N_n^{a_n}\smallsetminus W^u_n}$
defined by $N_I:=\{x\in N_n^{a_n}\mid x_n^s\in I\} $; set
$\la'^u_n:=\sup\{r'^u(\vp^u_n(x^u_n) \mid x\in N_I\}$ and 
$\la'^s_n:=\sup\{r'^s(\psi^s_n(x^s_n) \mid x\in N_I\}$.
If the product $\la'^u_n\la'^s_n$ is less than 1, the inequality $(*)_n$ is fulfilled  by the pair 
$(\varphi^u_n,\psi^s_n)$ and hence, the embedding $\phi_{\vp_n^u\psi_n^s}$ is well-defined on $N_n^{a_n}$.

If not, we replace $\psi^s_n$ with $f'^k\circ \psi^s_n$ with $k$ admissible and large enough. Indeed, 
the effect of this change is to multiply $\la'^s_n$ by some positive factor bounded
by  $(\frac 14)^k$ 
while $\la'^u_n$ is kept fixed and hence, 
$(*)_n$ becomes fulfilled when $k$ is large enough. Since the third item is empty for $i= n$, we have built 
some $\psi_n^s$ as desired.

For the induction, let us build $\psi^{s}_{i},~i<n,$  with the required properties 
assuming that the homeomorphisms  $\psi^{s}_{n},\dots,\psi^s_{i+1}$ have already been built. 
The stable manifolds of saddles in $\Si_{i}$ have heteroclinic intersections with unstable manifolds of saddles in $\Si_{j}$ with $j>i$ only. The image $\hat W^{s}_{i,i}$ of $W^s_i$ under $p_i:V_i\to \hat V_i$
is a  closed smooth 1-dimensional submanifold. According to Corollary \ref{con}, the intersection 
$\hat W^{s}_{i,i}\cap (\bigcup\limits_{j=i+1}^n\hat{N}^{a_{j}}_{j})$ consists of finitely many open arcs 
$\hat I^{i}_1,\dots,\hat I^{i}_{r_{i}}$ such that $\hat I^{i}_l$ for each $l=1,\dots,r_{i}$ is a connected 
component of $\hat W^{s}_{i,i}\cap\hat{N}^{a_{j}}_{j,i}$ for some $j>i$.  

In order to satisfy the third item of the statement, $\psi^s_i$ is defined on $p_{{i}}^{-1}(\hat I^{i}_l)$
in an equivariant way. Denote by $\psi^s_{i,l}$ this partial definition of $\psi^s_i$; its image is contained 
in $W'^s_{i,i}$. 

More precisely, if $I^i_{l,\al}$ is a connected component of $p_{{i}}^{-1}(\hat I^{i}_l)$
it is a proper arc in some $N^{a_j}_i$ and it intersects $W^u_j$ in a unique point $x^i_{l,\al}$.
Set $x'^i_{l,\al}= \vp^{us}(x^i_{l,\al})$ and denote $I'^i_{l,\al}$ the connected component 
of $W'^s_i\cap N'_j$ passing through the point $x'^i_{l,\al}$. Then, the restriction of $\psi^s_{i,l}$
to the arc $I^i_{l,\al}$ reads:
$$
\psi^s_{i,l,\al}= \phi_{\vp^u_j, \psi^s_j}\vert_{I^i_{l,\al}}: I^i_{l,\al}\to I'^i_{l,\al}\,.
$$

By Lemma \ref{pusk}, the map $\vp^{us}$ sends $W^{s}_i\cap L^u$ to $W'^{s}_i\cap {L'}^u$ preserving the order on each separatrix of $W^s_i\smallsetminus\Si_i$ and $W'^s_i\smallsetminus\Si'_i$. On the other hand, 
$\psi^s_{i,l,\al}$ is also order preserving. Both together imply that $\psi^s_{i,l}$ is order preserving 
since we know that it is an injective map. Moreover, the union of all $\psi^s_{i,l}$ -- which makes sense
as their respective domains are mutually disjoint -- is order preserving.
Therefore, there is an equivariant homeomorphism
$\psi^s_i: W^s_i\smallsetminus \Si_i\to  W^s_i\smallsetminus \Si_i$ which extends all $\psi^s_{i,l}$.

{Since $\vp^{us}$ is continuous, the above homeomorphism extends continuously to 
$\psi^s_i: W^s_i\to  W^s_i$.  At this point of the construction items 1 and 3 of the statement are satisfied.
The condition of item 2 follows from Lemma \ref{psibe} for  stable separatrices that contain heteroclinic 
points. If some stable separatrix has no heteroclinic points, one changes $\psi^s_i$ to $f'^k\circ \psi^s_i$ on the separatrix where $k$ is a large common multiple of the period of the separatrix, like to the construction made in the case $i=n$.}
\end{demo}\medskip

{\sc Proof of Theorem \ref{t.invariant} continued.}
Let us recall that we denoted by $\mathcal E:\mathbb R^3\to\mathbb R^3$ {the canonical  linear diffeomorphism}
with the unique fixed point {$O=(0,0,0)$} which is a saddle point 
whose  unstable manifold is the plane $Ox_1x_2$ 
 and stable manifold is the axis $Ox_3$; 
  for simplicity, we assume  that 
 $\mathcal E$ has a sign $\nu=+$ (see the beginning of Section \ref{section-compatible}). 
 Let $$N=\{(x_1,x_2,x_3)\in\mathbb R^3:0\leq(x_1^2+x_2^2)x_3\leq 1\}.$$  Let $\rho>0,\,\delta\in(0,\frac{\rho}{4})$ and $$d=\{(x_1,x_2,x_3)\in\mathbb R^3:x_1^2+x_2^2\leq\rho^2,x_3=0\},$$ $$U=\{(x_1,x_2,x_3)\in\mathbb R^3:(\rho-\delta)^2\leq x_1^2+x_2^2\leq\rho^2,x_3=0\},$$  $$c=\{(x_1,x_2,x_3)\in\mathbb R^3:x_1^2+x_2^2=\rho^2,x_3=0\},$$ $$c^{0}=\{(x_1,x_2,x_3)\in\mathbb R^3:x_1^2+x_2^2=(\rho-\frac{\delta}{2})^2,x_3=0\},$$ $$c^{1}=\{(x_1,x_2,x_3)\in\mathbb R^3:x_1^2+x_2^2=(\rho-\delta)^2,x_3=0\}.$$
Let $K=d\setminus int\,{\mathcal E}^{-1}(d)$, $V=(K\cup \mathcal E(K))\cap\{(x_1,x_2,x_3)\in\mathbb R^3:x_1\geq 0,x_3=0\}$ and $\beta=U\cap Ox^+_1$, where $Ox_1^+=\{(x_1,x_2,x_3)\in\mathbb R^3:x_3^2+x_2^2=0,x_1>0\}$. 

Choose a point $Z^0=(0,0,z^0)\in Ox_3^+$ such that $\rho^2 z^0<\frac{1}{4}$ (see Figure \ref{A}). Then,
choose a point $Z^1=(0,0,z^1)$ in $Ox_3^+$ so that $z^0>z^1>\frac{z^0}{4}$. Let  $\Pi(z)=\{(x_1,x_2,x_3)\in\mathbb R^3:x_3=z\}$.
In what follows,  for every subset $A\subset Ox_1x_2$, we denote by $\tilde A$ will denote the cylinder $\tilde A=A\times[0,z^0]$. Denote by {$\mathcal W$} the 3-ball bounded by the annulus $\tilde c$ and the two planes $\Pi(z^0)$ and $\Pi(\frac{z^{0}}{4})$. Let $\Delta$ be a closed 3-ball bounded by the surface $\tilde c^{1}$ and the two planes $Ox_1x_2$ and $\Pi(z^1)$. Let 
$$\mathcal T=\bigcup\limits_{k\in\mathbb Z}\mathcal E^{k}(\tilde d){\quad{\rm and}\quad}
\mathcal H=\bigcup\limits_{k\in\mathbb Z}\mathcal E^{k}(\Delta).$$
Notice that {the construction yields} $\mathcal H\subset int\,\mathcal T$ and makes {$\mathcal W$}
a fundamental domain for the action of $\mathcal E$ on $\mathcal T$. 

\begin{figure}[h]\centerline{\includegraphics[width=12cm,height=14cm]{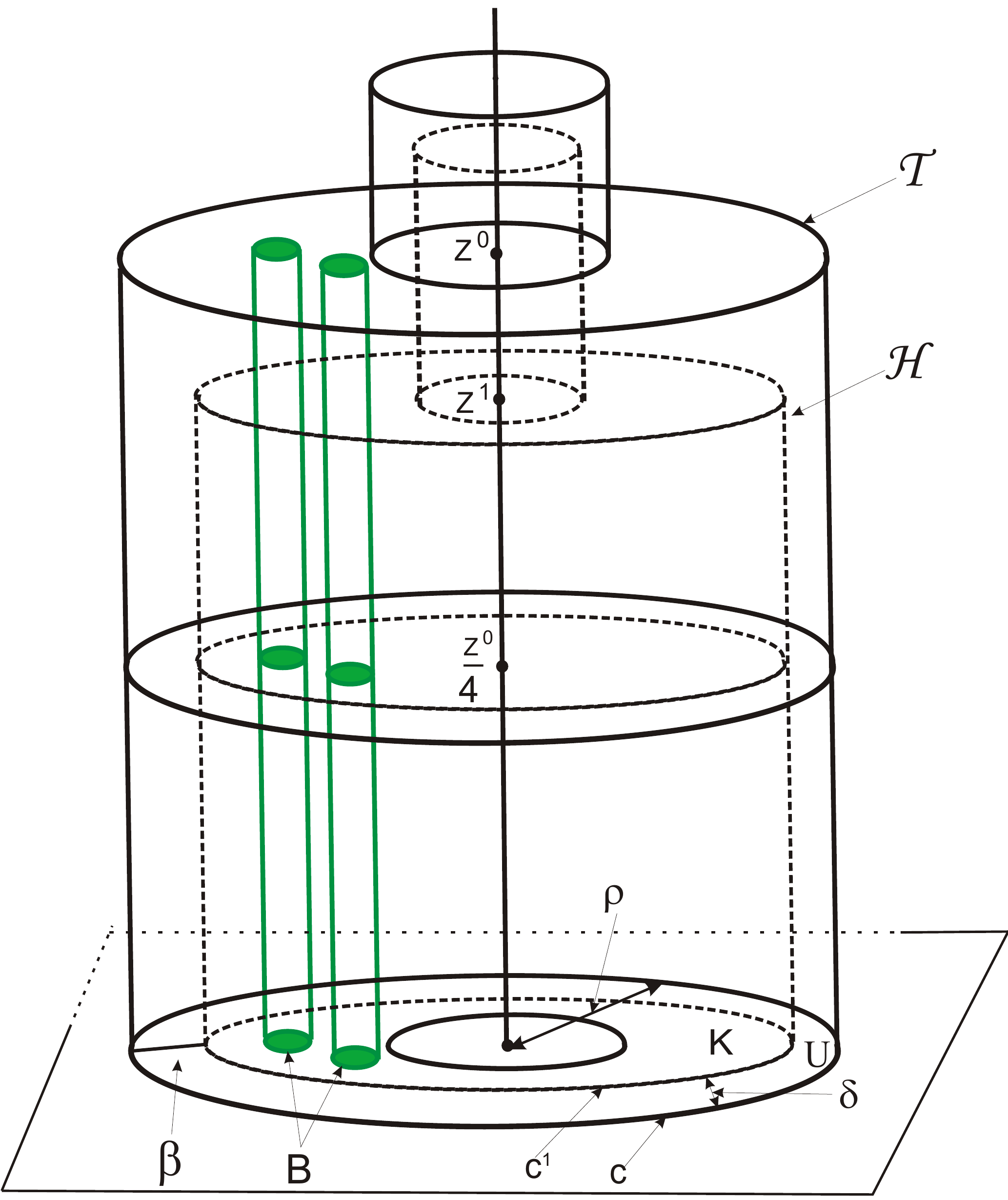}}\caption{A linear model} \label{A}
\end{figure}

Now, we come back to $f$ and construct some 
neighborhoods $\mathcal H_\ga\subset \mathcal T_\ga$ around each 
separatrix $\ga$ which contains heteroclinic points. Therefore,
we consider only the case $n\geq 1$ and separatrices of the saddle points from $\Si_i$ for $i\in\{0,\dots,n-1\}$, but 
not those from $\Si_n$ since their one-dimensional separatrices do not contain heteroclinic points.
Let $G_i$ be the union of all 
stable separatrices of saddle points {in} $\Si_i$ that contain heteroclinic points. Let $\check G_i\subset G_i$ be {the} union of separatrices {in} $G_i$ such that $G_i=\bigcup\limits_{\ga\in\check G_i}orb(\ga)$ and, for every pair $(\ga_1, \ga_2)$ of distinct separatrices in $\check G_i$ {and  every $k\in \mathbb Z$,} one has $\ga_{2}\neq f^k(\ga_1)$. {For $\ga\in G_i$ with the end point $p\in \Si_i$ and a point $q\in\Si_j,\,j>i$,
 let us consider a sequence of different periodic orbits $p=p_0\prec p_1\prec\dots\prec p_k=q$ such that $\ga\cap W^u_{p_1}\neq\emptyset$, the length of the longest such chain is denoted by $beh(q|\ga)$.}

Let $\ga\in\check G_i$ be a separatrix of $p\in \Si_i$ {and let $N^t_\ga$ be the connected component of $N^t_p \smallsetminus W^u_p$ which contains $\ga$}. We endow {with} the index $\ga$ {(resp. $p$)  the preimages in $M$ (through the linearizing map $\mu_p$) of  all objects from the linear model $\mathcal N$ associated with the separatrix $\ga$ (resp. $p$); for being precise we decide that $\mu_p(\ga)=Ox_3^+$.} For a separatrix $\ga$ in $\check G_i$, let us fix a saddle point $q_\ga$ such that $beh(q_\ga|\ga)=1$. {Notice that the intersection $\ga\cap W^u_{q_\ga}$ consists of a finite number {of} heteroclinic orbits.}

\begin{lemm}\label{order(n-1)} Let $n\geq 1$, {$i\in\{0,\dots,n-1\}$}. For every $\ga\in\check G_{i}$ there are  positive numbers $\rho$, $\de$ and $\ep$ (depending on $\ga$) such that for every heteroclinic point {$Z^0_\ga\in (\ga\cap W^u_{q_\ga})$} with $z^0<\ep$ the following properties hold:

{\rm (1)} $U_{p}$ avoids all heteroclinic points;

{\rm (2)} {$\vp(\tilde d_p)\subset {\phi}_{{\vp}^{u}_i,\psi^s_i}(N^{a_i}_i)$};

{\rm (3)} $\vp(\tilde c_p)\cap{\phi}_{{\vp}^{u}_i,\psi^s_i}(\tilde c^0_p)=\emptyset$, $\vp(\tilde c_p^1)\cap{\phi}_{{\vp}^{u}_i,\psi^s_i}(\tilde c_p^0)=\emptyset$ and $\vp(\tilde\beta_\ga)\subset{\phi}_{{\vp}^{u}_i,\psi^s_i}(\tilde V_\ga)$.
\end{lemm}    
\begin{demo} Let $\ga\in{\check G_{i}},\,i\in\{0,\dots,n-1\}$. Due to Lemma \ref{Ge}, there is a 
generic  $\rho>0$ such that the  curve $c_\ga$ avoids all heteroclinic points. Since ${W^s_{l}}$ accumulates on $W^s_{k}$ for every $l<k$, then $K_{p}\cap W^s_{i-1}$ is made of a finite number of heteroclinic points $y_1,\dots,y_r$ which we can cover {by} closed 2-discs $b_1,\dots,b_r\subset int\,K_p$. In $K_{p}\smallsetminus int(b_1\cup\dots\cup b_r)$ there is a finite number of heteroclinic points from $W^s_{i-2}$ which we cover by the union of a finite number {of} closed 2-discs, and so on. Thus we get that all heteroclinic points in $K_{p}$ belong to the union of {finitely many} closed 2-discs avoiding $\partial K_{p}$. Therefore, there is $\delta\in(0,\frac{\rho}{4})$ such that  $U_{p}$ avoids heteroclinic points. This proves item (1). 

By assumption of Theorem \ref{t.invariant}, $\vp$ is defined on the complement 
of the stable manifolds and, by definition, ${\phi}_{{\vp}^{u}_{i},\psi^s_{i}}$ coincides with $\vp$ on $W^u_i\smallsetminus L^s$, and hence on $U_p$. As $\vp$ and ${\phi}_{{\vp}^{u}_{i},\psi^s_{i}}$ are continuous, we can choose $\ep>0$ sufficiently  small so that, if $Z^0_{\ga}$ is any  heteroclinic point in the intersection $\ga\cap W^u_{q_\ga}$ with $z^0<\ep$, {the requirements of (2) and (3) are fulfilled.}
\end{demo}

Let us fix $U_p$ satisfying  item (1) of Lemma \ref{order(n-1)} and define
$$U_i=\bigcup\limits_{p\in \Si_i}\left(\bigcup\limits_{k=0}^{per(p)-1}f^k(U_{p})\right),\quad K_i=\bigcup\limits_{p\in \Si_i}\left(\bigcup\limits_{k=0}^{per(p)-1}f^k(K_{p})\right).$$
Until the end of Section \ref{IV}, we assume that  for every $\ga\in \check G_{n-1}$ the neighborhoods $\mathcal T_\ga$ and $\mathcal H_\ga$ have  the parameters $\rho, \de,\ep, z^0$ as in Lemmas \ref{order(n-1)} and $z^1$ is chosen such that the arc $(z^0_\ga,z^1_\ga)\subset\ga$ does not contain heteroclinic points. But, when $\ga\in \check G_{i},\  i\in\{0,\dots,n-2\}$,  the parameter $\ep$ will be still more specified in Lemma \ref{order} below.

\begin{lemm}\label{order} Let $n\geq 2$. For every $i\in\{0,\dots,{n-2}\}$ and $\ga\in\check G_i$, there is a heteroclinic point $Z^0_\ga\in \ga$ satisfying the conditions of Lemma \ref{order(n-1)} and in addition:
  $$\mathcal T_\ga\cap\tilde U_j=\emptyset~~~for~~~j\in\{i+1,\dots,{n-1}\}.$$
\end{lemm}  

In this statement, it is meant that  $\tilde U_{n-1}$ is associated with the points $Z^0_{\ga'}, \ga'\in \check G_{n-1}$ given by Lemma \ref{order(n-1)} and $\tilde U_{j}$ is associated with the points 
$Z^0_{\ga''}, \ga''\in\check G_{j}$ given by Lemma \ref{order} {for every $j>i$}. 
Therefore, it makes sense to prove Lemma \ref{order} by decreasing induction on $i$  from ${i=n-2}$ to $0$. This is what is done below.\\ 

\nd\begin{demo} Let us first prove the lemma  for $i= n-2$. Let $\ga\in {\check G_{n-2}}$ and let $p$ be the saddle end point of $\ga$. Notice that the intersection $\ga\cap K_{n-1}$ consists of a finite number points $a_1,\dots,a_l$ avoiding $U_{n-1}$. Let $d_{1},\dots,d_{l}\subset K_{n-1}$ be compact discs with centres $a_1,\dots,a_l$ and radius $r_*$ (in linear coordinates of $N_p$) avoiding $U_{n-1}$. Let us choose a number $n^*\in\mathbb N$ such that $\frac{\rho}{2^{n^*}}<r^*$. Let $Z^*_\ga\subset\ga$ be a point such that the segment $[p,Z^*_\ga]$ of $\ga$ avoids $\tilde K_{n-1}$ and $\mu_p(Z^*_\ga)=Z^*=(0,0,z^*)$ where $z^*<\ep$. Then every heteroclinic point $z^0_{\ga}$ so that $z^0<\frac{z^*}{2^{n^*}}$ possesses the property: $\mathcal T_\ga\cap\tilde K_{n-1}$ avoids $\tilde U_{n-1}$. 

 {For the induction, let us assume now that the construction of the desired heteroclinic points}
is done for $i+1,i+2,\dots,n-2$.
 Let us do it for $i$. 
Let $\ga\in  {\check G_i}$. By assumption of the induction $(\bigcup\limits_{k=i+1}^{j-1}\mathcal T_k)\cap\tilde U_j=\emptyset$ for $j\in\{i+2,\dots,n-1\}$. Since ${W^s_{k-1}}$ accumulates on $W^s_{k}$ for every $k\in\{0,\dots,n\}$,
 then $(\bigcup\limits_{k=i+1}^{j-1}\mathcal T_k)\cap K_j$ is a compact subset of $K_j$ and the intersection $(\ga\setminus(\bigcup\limits_{k=i+1}^{j-1}\mathcal T_k))\cap K_j$ consists of a finite number points $a_1,\dots,a_l$ avoiding $U_{j}$. Let $d_{1},\dots,d_{l}\subset K_{j}$ be compact discs with centres $a_1,\dots,a_l$ and radius $r_*$ (in linear coordinates of $N_p$) avoiding $U_{j}$ and such that $r_*$ is less than  {the distance between} $\partial(K_j\setminus U_j)$ and $(\bigcup\limits_{k=i+1}^{j-1}\mathcal T_k)\cap K_j$. Similar to the case $i=n-2$ it is possible to choose a heteroclinic point $Z^0_{\ga}$ sufficiently close to the saddle $p$ where $\ga$ ends such that the set $(\mathcal T_\ga\setminus(\bigcup\limits_{k=i+1}^{j-1}\mathcal T_k))\cap \tilde K_j$ avoids $\tilde U_j$. 
\end{demo}
 {\medskip}

In what follows, we assume that, for every $\ga\subset\check G_i,\,i\in\{0,\dots,n-2\}$, the neighborhoods  
$\mathcal T_\ga$ and $\mathcal H_\ga$ have parameters $\rho, \de,\ep, z^0$ as in Lemma \ref{order(n-1)} and moreove $\ep$ satisfies to Lemma \ref{order}.
For $i\in\{0,\dots,n-1\}$, we set
$$\mathcal T_i=\bigcup\limits_{\ga\subset\check G_i}\left(\bigcup\limits_{k=0}^{per(\ga)-1}f^k(\mathcal T_{\ga})\right).$$ 

For $\ga\subset\check G_i,\, j>i$, let us denote by ${\mathcal{J}}_{\ga,j}$ the union of all connected components of ${W}^u_{j}\cap\mathcal T_\ga$ which do not lie in {$int\,\mathcal T_k$} with $i<k<j$. 
Let $\mathcal{J}_\ga=\bigcup\limits_{j=i+1}^{n}{\mathcal{J}}_{\ga,j}$ {and $$\mathcal J_i=\bigcup\limits_{\ga\subset\check G_i}\left(\bigcup\limits_{k=0}^{per(\ga)-1}f^k(\mathcal J_{\ga})\right).$$} 

Let {$\mathcal W_\ga$} be the fundamental domain of $f^{per(\ga)}\vert_{\mathcal T_\ga\setminus W^u_p}$
 {limited by the plaques of the two heteroclinic  points $Z^0_\ga$ and $f^{per(\ga)}Z^{0}_\ga$. Notice that}
 $\ga\cap{\mathcal W_\ga}$ is a fundamental domain of $f^{per(\ga)}\vert_\ga$. Since ${W^u_{k}}$ accumulates on $W^u_{l}$ {only when} 
 $l<k$, then the set ${\mathcal{J}}_{\ga,j}\cap{\mathcal W_\ga}$ consists of a finite number  {of}  closed 2-discs. Hence, the set $\mathcal{J}_\ga\cap\ga\cap{\mathcal W_\ga}$ consists of a finite number 
  {of heteroclinic points;} denote them $Z^2_\ga, \ldots, Z^{m}_\ga$ ($m$ depends on $\ga$).  {Finally, choose an} arbitrary point $Z^1_\ga\in\ga$ so that the arc $(z^0_\ga,z^1_\ga)\subset\ga$ does not contain heteroclinic points from $\mathcal J_\ga$. 
 Let us construct $\mathcal H_\ga$ using the point $Z^1=\mu_p(Z^1_\ga)$ and the parameter $\de$ from Lemma \ref{order(n-1)}.
 Without loss of generality we will assume that $\mu_p(Z^i_\ga)=Z^i=(0,0,z^i)$ for $z^0>z^1>\dots>z^m>\frac{z^0}{4}$. For $i=0, \ldots, n-1$ let $$\mathcal H_i=\bigcup\limits_{\ga\subset\check G_i}\left(\bigcup\limits_{k=0}^{per(\ga)-1}f^k(\mathcal H_{\ga})\right){\quad{\rm and}\quad}\mathcal {M}_{i}=V_f\cup\bigcup\limits_{k=0}^{i}{(G_k\cup \Si_k)}.$$

{\begin{lemm} \label{main0} There is an equivariant topological embedding ${\varphi}_0:\mathcal{M}_0\to {M}'$ with following properties:

{\rm (1)} ${\varphi}_0$ coincides with ${\varphi}$ out of $\mathcal{T}_0$;

{\rm (2)} ${\varphi}_0\vert_{\mathcal H_0}={\phi}_{\psi^u _{0},\psi^s_{0}}\vert_{\mathcal H_0}$,
 where $\psi^u_0={\varphi}\vert_{W^u_0}$;

{\rm (3)} $\vp_0(W^u_1)=W'^u_1$ and $\vp_0(W^u_k\setminus\bigcup\limits_{j=1}^{k-1}int\,\mathcal T_j)\subset W'^u_k$ for every $k\in\{2,\dots,n\}$.
\end{lemm}}
\nd\begin{demo} The desired $\vp_0$ should be an interpolation between 
 $\vp: V_f\smallsetminus \mathcal T_0 \to  {M'}$ and $\phi_{{\vp}^u _{0},\psi^s_{0}}\vert_{\mathcal H_0}$.
  Due to Lemma \ref{order(n-1)}  (2) {and the equivariance of the considered maps}, the embedding 
$${\xi_0}={\phi}^{-1}_{\psi^{u}_0,\psi^s_0} {\vp}:\mathcal T_0\setminus W^s_0\to M$$ is well-defined. 
Let $\ga\subset\check G_0$ be a separatrix ending at $p\in\Si_0$ and $\xi_\ga=\xi_0|_{\mathcal T_\ga}$. 
 {By construction,} the topological embedding $\xi=\mu_p\xi_\ga\mu^{-1}_p:\mathcal T\to N$ has the following properties:

(i) $\xi\mathcal E=\mathcal E\xi$ (as $\mathcal E\mu_p=\mu_pf^{per(\ga)}$ and $\xi_\ga f^{per(\ga)}=f^{per(\ga)}\xi_\ga$);

(ii) $\xi$ is the identity on $Ox_1x_2$ (as ${\phi}_{\psi^{u}_0,\psi^s_0}|_{W^u_0}=\vp|_{W^u_0}$);

(iii) $\xi(\Pi(z^0)\cap\mathcal T)\subset \Pi(z^0)$ and $\xi(\Pi(z^i)\cap\partial\mathcal T)\subset \Pi(z^i),\,i\in\{2,\dots,m\}$ (as {$\xi_\ga(L^u\cap\mathcal T_\ga\setminus\ga)\subset L^u$});

(iv) $\xi(\tilde c)\cap\tilde c^0=\emptyset$, $\xi(\tilde c^1)\cap\tilde c^0=\emptyset$ and $\xi(\tilde\beta)\subset\tilde V$ (due to Lemma \ref{order(n-1)} (3)). 

 \nd Thus, $\xi$ satisfies all conditions of Proposition \ref{F1} and, hence, there is an embedding $\zeta:\mathcal T\to N$ such that:   

(I) $\zeta a=a\zeta$;

(II) $\zeta$ is the identity on $\mathcal H$;

(III) $\zeta(\Pi(z^i)\cap\mathcal T)\subset \Pi(z^i),\,i\in\{0,2,\dots,m\}$

(IV) $\zeta$ is $\xi$ on $\partial{\mathcal T}$.

\nd Then the embedding $\zeta_\ga=\mu_p^{-1}\zeta\mu_p:\mathcal T_\ga\to N_\ga$ satisfies the properties:

(I') $\zeta_\ga f^{per(\ga)}=f^{per(\ga)}\zeta_\ga$;

(II') $\zeta_\ga$ is the identity on $\mathcal H_\ga$;

(III') $\zeta( {\mathcal J_\ga)}\subset L^u$

(IV') $\zeta_\ga$ is $\xi_\ga$ on $\partial{\mathcal T}_\ga$.

\nd Independently, one does the same for every separatrix $\ga\subset\check G_0$ {. Then, it is extend
to all separatrices in $G_0$ by equivariance.} {As a result}, we get a homeomorphism $\zeta_0$  {of $\mathcal T_0$ onto
$\xi_0(\mathcal T_0)$ which coincides with $\xi_0$ on $\partial\mathcal T_0$. Now, define the embedding 
$\vp_0:\mathcal M_0\to M'$ to be equal to ${\phi}_{\psi^{u}_0,\psi^s_0}\zeta_0$ on $\mathcal T_0$ 
and to $\vp$ on $\mathcal M_0\setminus \mathcal T_0$. One checks the next properties:}

$(1)\, {\varphi}_0$ coincides with ${\varphi}$ out of $\mathcal{T}_0$;

$(2)\, {\varphi}_0\vert_{\mathcal H_0}={\phi}_{\psi^u _{0},\psi^s_{0}}\vert_{\mathcal H_0}$;

$(3')\, \vp_0( {\mathcal J}_0)\subset L^u$.

The last property and the definition of the set $ {\mathcal J}_\ga$ imply that $\vp_0(W^u_1)=W'^u_1$ and $\vp_0(W^u_k\setminus\bigcup\limits_{j=1}^{k-1}int\,\mathcal T_j)\subset W'^u_k$ for every $k\in\{2,\dots,n\}$.
Thus $\vp_0$ satisfies  {all  required conditions of the} lemma. 
\end{demo}

\begin{lemm} \label{main} Assume $n\geq 2$, $i\in\{0,\dots,n-2\}${, and assume} there is an equivariant topological embedding ${\varphi}_i:\mathcal{M}_i\to{M'}$ with following properties:

{\rm (1)} ${\varphi}_i$ coincides with ${\varphi}_{i-1}$ out of $\mathcal{T}_i$;

{\rm (2)} ${\varphi}_i\vert_{\mathcal H_i}={\phi}_{\psi^u _{i},\psi^s_{i}}$,
 where $\psi^u_i={\varphi}_{i-1}\vert_{W^u_i}$ {and} $\vp_{-1}=\vp$; 

{\rm (3)} there is an $f$-invariant {union of tubes}
 $\mathcal B_i\subset(\mathcal T_i\cap\bigcup\limits_{j=0}^{i-1}\mathcal H_j)$ {containing
 $(\mathcal T_{i}\cap(\bigcup\limits_{j=0}^{i-1}W^s_j))$} where ${\varphi}_i$ coincides with $\vp_{i-1}$ 
 $($we assume $\mathcal B_0=\emptyset)$;
 
{\rm (4)} $\vp_i(W^u_{i+1})=W'^u_{i+1}$ and $\vp_i(W^u_k\setminus\bigcup\limits_{j=i+1}^{k-1}int\,\mathcal T_j)\subset W'^u_k$ for every $k\in\{i+2,\dots,n\}$.

{\nd Then} there is a homeomorphism $\vp_{i+1}$ with the same properties {\rm (1)-(4)} 
\end{lemm}
\nd \begin{demo} The desired $\vp_{i+1}$ should be an interpolation between 
 $\vp_i: \mathcal M_{i+1}\smallsetminus \mathcal T_{i+1} \to M'$ and 
 $\phi_{{\psi}^u_{i+1},\psi^s_{i+1}}\vert_{\mathcal H_{i+1}}$ {where $\psi_{i+1}^u=\vp_i\vert_{W^u_{i+1}}$.}
Let $\ga\subset\check G_{i+1}$ be a separatrix ending at $p\in\Si_{i+1}$. {It follows from the definition of the set $\mathcal J_i$ and the choice of the point $q_\ga$ that $(W^u_{q_\ga}\cap \mathcal T_i)\subset \mathcal J_i$. Then, due to condition (4) for $\vp_i$ we have $\vp_i(W^u_{q_\ga}\cap \mathcal T_i)\subset W^u_{q'}$. By the property (1) of the homeomorphism $\vp_i$ and the properties of $\mathcal T_{i+1}$ from Lemmas \ref{order(n-1)} {(1)} and \ref{order}, we get that $\vp_i|_{\tilde U_p}=\vp|_{\tilde U_p}$. Then ${\phi}_{{\vp}^{u}_{i+1},\psi^s_{i+1}}|_{\tilde U_p}={\phi}_{{\psi}^{u}_{i+1},\psi^s_{i+1}}|_{\tilde U_p}$. Thus it follows from the property (2) in Lemma \ref{order(n-1)} that the following embedding is well-defined: $\xi_\ga={\phi}_{{\psi}^{u}_{i+1},\psi^s_{i+1}}^{-1}\vp_i:{\mathcal T}_\ga\setminus(\ga\cup p)\to M'$.}

{By construction,} the topological embedding $\xi=\mu_p\xi_\ga\mu^{-1}_p$ satisfies to all conditions of Proposition \ref{F1}. Let $\zeta$ be the  embedding {which is yielded by that proposition}.
Define $\zeta_\ga=\mu_p^{-1}\zeta\mu_p$. Notice that by the property (3) of the homeomorphism $\psi^s$ in Lemma \ref{1dim} and by the properties $\psi^u_{i+1}=\vp_i|_{W^u_i}$, we have that $\zeta_{\ga}$ is {the} identity on a neighborhood $\tilde B_{\ga}\subset(\mathcal T_\ga\cap\bigcup\limits_{j=0}^i\mathcal H_j)$ of $\mathcal T_{\ga}\cap(\bigcup\limits_{j=0}^iW^s_j)$. Independently, one does the same for every separatrix $\ga\subset\check G_{i+1}$. {Assuming that $\zeta_{f(\ga)}=f'\zeta_{\ga}f^{-1}$  and $\tilde B_{i+1}=\bigcup\limits_{\ga\subset\check G_{i+1}}\left(\bigcup\limits_{k=0}^{per(\ga)-1}f^k(\tilde B_{\ga})\right)$ we get a homeomorphism $\zeta_{i+1}$ on $\mathcal T_{i+1}$.}  Thus the required homeomorphism coincides with ${\phi}_{\psi^{u}_{i+1},\psi^s_{i+1}}$ on $\mathcal H_{i+1}$ and with $\vp_i$ out of $\mathcal T_{i+1}$. 
\end{demo}

{Let $G$ be the union of all stable one-dimensional separatrices which do not contain heteroclinic points, $ N^t_{G}=\bigcup\limits_{\gamma\subset G} N^t_\gamma$ and  $$\mathcal {M}=\mathcal M_{n-1}\cup G.$$ 

{\begin{lemm} \label{main-} There are {numbers} 
$0<\tau_1<\tau_2<1$ and an equivariant embedding $h_{\mathcal M}:\mathcal M\to {M'}$ with the following properties: 

{\rm (1)} $h_{\mathcal M}$ coincides with {${\varphi}_{n-1}$} 
out of $N^{\tau_2}_{G}$;

{\rm (2)} $h_{\mathcal M}$ coincides with ${\phi}_{{\vp}_{n-1},{\psi}^s}$ on $\vert_{\mathcal N^{\tau_1}_{G}}$,
{where $\psi^s: L^s\to L'^s$ is yielded by 
Lemma \ref{1dim};}

{\rm (3)} there is an $f$-invariant neighborhood of the set {$N_{G}\cap(G_0\cup\dots\cup G_{n-1})$} where $h_{\mathcal M}$ coincides with {$\vp_{n-1}$.}
\end{lemm}}
\begin{demo} Let $\check G\subset G$ be a union of separatrices from $G$ such that $\ga_{2}\neq f^k(\ga_1)$ for every $\ga_1,\ga_2\subset\check G$, $k\in\mathbb Z\setminus\{0\}$ and $G=\bigcup\limits_{\ga\in\check G}orb(\ga)$. Let $i\in\{0,\dots,n\}$, $p\in\Si_i$ and $\gamma\subset G$. 

Notice that $\left(N_{\ga}\setminus(\ga\cup p)\right)/f^{per(\ga)}$ is homeomorphic to $X\times[0,1]$ where 
$X$ is 2-torus and the natural projection $\pi_{\ga }: N_{\ga }\setminus(\ga \cup p)\to X\times[0,1]$ sends 
$\partial N^t_{\ga }$ to $X\times\{t\}$ for each $t\in(0,1)$ and sends $W^u_p \setminus p$ to $X\times\{0\}$. 
Let  $\xi_{\ga }={\phi}^{-1}_{{\vp}_{n-1}|_{W^u_i},{\psi}^s_i}\vp_{n-1}|_{N_{\ga }^{a_i}\setminus(\ga \cup p)}$ and  $\hat\xi_{\ga}=\pi_{\ga }\xi_{\ga }\pi^{-1}_{\ga }|_{X\times[0,a_i]}$. Due to item (3) of Lemma \ref{main}, the homeomorphism $\hat\xi_{\ga }$ coincides with {the} identity in some neighborhood of 
$\pi_\ga(N^{a_i}_{\ga}\cap(G_0\cup\dots\cup G_{n-1}) )${. Let} us choose {this neighborhood} of the form $B_\ga\times[0,a_i]$. Let us choose numbers $0<\tau_{1,\ga }<\tau_{2,\gamma }<a_i$ such that 
$\hat \xi_{\ga }(X\times[0,{\tau_{2,{\ga }}}])\subset X\times[0,\tau_{1,\ga })$. 
{By construction,} $\hat \xi_{\ga }:X\times[0,{\tau_{2,{\ga }}}]\to {X\times[0,1]}$ {is} a topological embedding which is the identity on $X\times\{0\}$ and  
$\hat \xi_{\ga }|_{B_{\ga }\times[0,\tau_{2,\ga }]}=id|_{B_{\ga }\times[0,\tau_{2,\ga }]}$. 
Then, by 
 Proposition \ref{F2}, 

1. there is a homeomorphism $\hat \zeta_{\ga }:X\times[0,\tau_{2,\ga }]\to \hat\xi(X\times[0,\tau_{2,\ga }])$ such that  $\hat\zeta_{\ga }$ is identity on $X\times[0,{\tau_{1,{\ga }}}]$ and is $\hat \xi_{\ga }$ on $X\times\{\tau_{2,\ga}\}$. 

2. $\hat \zeta_{\ga }|_{B_{\ga }\times[0,\tau_{2,\ga }]}=id|_{B_{\ga }\times[0,\tau_{2,\ga }]}$.

Let $\zeta_{\ga }$ be a lift of $\hat \zeta_{\ga }$ on $N^{\tau_{2,\ga}}_{\ga }$ which $\xi_\ga$ on $\partial N^{\tau_{2,\ga}}_{\ga }$. Thus $\vp_{\ga }={\phi}_{{\vp}_{n-1}|_{W^u_i},{\psi}^s_i}\zeta_{\ga }$ is the desired extension of $\vp_{n-1}$ to ${ N}_{\ga }$. Doing the same for every separatrix $\ga \subset \check G$ and {extending it to the other separatrices from $G$ by equivariance,}
we get the {required} homeomorphism $h_{\mathcal M}$ for $\tau_1=\min\limits_{\gamma \subset\check G}\{\tau_{1,\ga }\}$ and $\tau_2=\min\limits_{\gamma \subset\check G}\{\tau_{2,\ga }\}$.
\end{demo}

So far in this section, we have modified the homeomorphism $\varphi$ in a union of suitable linearizable
neighborhoods of $\Om_2$ (with  their 1-dimensional separatrices removed) so that the modified homeomorphism
extends equivariantly to $W^s(\Om_2)$. At the same time, we can perform  a similar modification about 
$\Om_1$ since the involved linearizable neighborhoods of $\Om_2$ and $\Om_1$ are mutually disjoint. {Thus,} we get a homeomorphism $h:M\setminus(\Om_0\cup\Om_3)\to M\setminus(\Om'_0\cup\Om'_3)$ conjugating $f|_{M\setminus(\Om_0\cup\Om_3)}$ with $f'|_{M\setminus(\Om'_0\cup\Om'_3)}$. Notice that $M\setminus(W^s_{\Om_1}\cup W^s_{\Om_2}\cup\Om_3)=W^s_{\Om_0}$ and $M\setminus(W^s_{\Om'_1}\cup W^s_{\Om'_2}\cup\Om'_3)=W^s_{\Om'_0}$. Since $h(W^s_{\Om_1})=W^s_{\Om'_1}$ and 
$h(W^s_{\Om_2})=W^s_{\Om'_2}$, then $h(W^s_{\Om_0}\setminus\Om_0)=W^s_{\Om'_0}\setminus\Om'_0$. Thus for each connected component $A$ of $W^s_{\Om_0}\setminus\Om_0$, there is a sink $\omega\in\Om_0$ such that $A=W^s_\omega\setminus\omega$. Similarly, $h(A)$ is a connected component of $W^s_{\Om'_0}\setminus\Om'_0$ such that  $h(A)=W^s_{\omega'}\setminus\omega'$ for a sink $\omega'\in\Om'_0$. Then we can continuously extend $h$ to $\Om_0$ {by defining} $h(\omega)=\omega'$ for every $\omega\in\Om_0$. A similar extension of $h$ to $\Om_3$ finishes the proof {of Theorem \ref{t.invariant}}.

\section{Topological background}\label{V}
{We use below the notations which have been} introduced before Lemma \ref{order(n-1)}.

\begin{prop} \label{F1} Let $z^0>z^1>\dots>z^m>\frac{z^0}{4}>0$ and $\xi:\mathcal T\setminus Ox_3\to N$ be a topological embedding with the following properties:

{\rm (i)} $\xi\mathcal E=\mathcal E\xi$;

{\rm (ii)} $\xi$ is the identity on $Ox_1x_2$;

{\rm (iii)} $\xi(\Pi(z^0\cap \mathcal T)) =\Pi(z^0)$ and $\xi(\Pi(z^i)\cap\partial\mathcal T)\subset \Pi(z^i),\,i\in\{2,\dots,m\}$;

{\rm (iv)} $\xi(\tilde c)\cap\tilde c^0=\emptyset$, $\xi(\tilde c^1)\cap\tilde c^0=\emptyset$ and $\xi(\tilde\beta)\subset\tilde V$. 

\nd Then there is a homeomorphism $\zeta:\mathcal T\to N$ such that  

{\rm (I)} $\zeta\mathcal E=\mathcal E\zeta$;

{\rm (II)} $\zeta$ is the identity on $\mathcal H$ {-- and consequentively on $Ox_1x_2$};

{\rm (III)} $\zeta(\Pi(z^i)\cap\mathcal T)\subset \Pi(z^i),\,i\in\{0,2,\dots,m\}$

{\rm (IV)} $\zeta$ is $\xi$ on $\partial{\mathcal T}$.

\nd Moreover, if $\xi$ is identity on $\tilde B$ for a set 
$B\subset (K\setminus U)$ then $\zeta$ is also identity on $\tilde B$.
\end{prop} 
\begin{demo} For $j= 0, ..., m$, we denote by $E_j$ the domain of $\R^3$ located between the horizontal
planes $\Pi(z_j)$ and $\Pi(z_{j+1})$, with $z_{m+1}= \frac {z_0}{4}$. Since the requested $\zeta$ has to be 
equivariant with respect to $\mathcal E$, it is useful to describe a fundamental domain $\mathcal V$ for the action of $\mathcal E$ on the closure of $\mathcal T\setminus \mathcal H$; the natural one is 
 $$\mathcal V=cl( \mathcal T\setminus \mathcal H)\cap(\bigcup_{j=0}^m E_j),$$
where $cl(-)$ stands for {\it closure of $(-)$}. The domain $\mathcal V$ is sliced by the horizontal
planes $\Pi(z_j), j= 2,\ldots, m,$ and the vertical cylinders $\mathcal E^{-1}(\tilde c)$ and $\tilde c^1$,
yielding the decomposition $\mathcal V= R_0\cup R_1\cup Q_0\cup Q_2\cup \ldots \cup Q_m$  
into solid tori whose interiors are mutually disjoint. Notice that the plane $\Pi(z_1)$ is not used in this decomposition.

More precisely, $R_0\subset E_0$ is limited by the cylinders $\mathcal E^{-1}(\tilde c^1)$ and $\mathcal E^{-1}(\tilde c)$; then, $R_1\subset E_0$  is limited by the cylinders $\mathcal E^{-1}(\tilde c)$ and $\tilde c^1$. The others of the list
are obtained from $\tilde U$ by slicing $\mathcal V$
with horizontal planes. The first of the latter, namely  $Q_0$, is special as 
it is bounded by $\Pi(z_0)$ and $\Pi(z_2)$;  then, $Q_j$ is bounded by $\Pi(z_j)$ and $\Pi(z_{j+1})$
for $j=2,...,m$. The vertical parts in the boundaries of the above-mentioned solid tori are provided  by
the vertical slices or the vertical parts of $\partial\mathcal T\cup\partial\mathcal H $. 

For $j= 0, 2, \ldots, m$, let $U(z_j):= \tilde U\cap\Pi(z_j)$. 
By construction, the top face of $R_0$ is\break 
$U'(z_0):=\mathcal E^{-1}\left(U(z_{m+1})\right)=  \Pi(z_0)\cap\mathcal E^{-1}(\tilde U)$; its bottom  is 
$U'(z_1):= \Pi(z_1)\cap\mathcal E^{-1}(\tilde U)$. Similarly, the top of $R_1$ is 
$U''(z_0):=  \Pi(z_0)\cap\tilde K $ and its bottom is $U''(z_1):=  \Pi(z_1) \cap\tilde K $.

It is important that each horizontal or vertical annulus  $Ann$ from the
previous list is marked with a special arc noted $\beta(Ann)$ linking the two boundary 
components of $Ann$ and defined as follows:

$$\beta(Ann)= Ann\cap \{x_1>0, x_2=0\}.$$ 
According to assumption (iv), all these arcs (except when $Ann= U''(z_0)\  {\rm or}\ U''(z_1)$) 
fulfill the next condition, referred to as the $\beta$-condition, namely: they are 
mapped by $\xi$ into $\{x_1>0\}$.

First of all, we define $\zeta\vert_{R_1}$ by rescaling $\zeta\vert_{\tilde K}$ in the next way. There is a
homeomorphism $\kappa: \tilde K\to R_1$ of the form: $(x_1,x_2,x_3)\mapsto (x_1,x_2, \rho(x_3))$
where $\rho: [0, z_0]\to [z_1,z_0]$ is any increasing continuous function. Then, we 
define $\zeta\vert_{R_1}= \kappa\,\xi\vert_{\tilde K}\,\kappa^{-1}$. Observe that $\zeta$
equals $\xi$ on $U''(z_0)$ and coincide with the identity on $U''(z_1)$. As a consequence,
the complement part of the statement follows directly. Indeed, if $B$ lies in $K$ and $\xi\vert_{\tilde B}= Id$
then its conjugate by $\kappa$ is the identity on $\tilde B\cap R_1$.

We continue by defining $\zeta$ on the other horizontal annuli from the  previous list. As required, $\zeta$
is the identity when this annulus lies in $\mathcal H$. 
For the others, that is, $U'(z_0)$ and $U(z_j),\ j= 0,2,...,m$, Lemma \ref{2-annulus} is applicable as it is 
explained right below.

Each of these annuli is bounded by two curves; one of the two lies in the frontier of $\mathcal T$
on which $\zeta$ has to coincide with $\xi$ and is mapped in the respective plane $\Pi(z_j)$
 -- according to (iii); and the other 
 lies in $\mathcal H$ where $\zeta$ has to coincide with $Id\vert_{\mathcal H}$. 
 In order to satisfy the equivariance property 3), we choose 
 $$(*)\quad\quad\zeta\vert_{U(z_{m+1})}=\mathcal E\,\zeta\vert_{U'(z_{0})}\,\mathcal E^{-1}.$$
 Moreover, due to  the {\it $\beta$-condition}, 
Lemma \ref{2-annulus} holds  and yields $\zeta$ on each of the listed horizontal 
annuli. 

We continue by defining $\zeta$ on the vertical annuli in the above splitting of $\mathcal V$.
When such an annulus lies in $\partial H$ (resp. $\partial\mathcal T$), we must take $\zeta=Id$ 
(resp. $\zeta=\xi$) over there.
The last two annuli are $R_0\cap R_1$ and $R_1\cap Q_0$ on which $\zeta$ is already defined 
by conjugating by $\kappa$. Notice that the $\beta$-condition holds for these two annuli 
because conjugating
by $\kappa$ preserves the $\beta$-condition. 

Let us look more precisely to $\partial Q_0$. It is made of the following: two horizontal annuli
$U(z_0)$ and $U(z_2)$, and three vertical ones $R_1\cap Q_0$, $\tilde c^1\cap E_1$
(lying in $\mathcal H$) and $\tilde c\cap (E_0\cup E_1)$. The $\beta$-condition holds for each of these
latter annuli.

For finishing the proof, it remains to extend the $\zeta$ which we have defined right above on the tori
$\partial R_0$, $\partial Q_0$ and 
 $\partial Q_j,\ j= 2, ..., m$ to embeddings of the solid tori $R_0, Q_0, Q_2,...,Q_m$
  with values in $E_0$, $E_0\cup E_1$, $ E_2,...,E_m$ respectively. The boundary data consists of annuli 
 where $\zeta $ fulfills the $\beta$-condition. 
 Therefore, the assumptions of Proposition \ref{solid} 
 are fulfilled; then the conclusion holds true and yields the desired extention of $\zeta$ to the listed solid tori.
 
 According to $(*)$, we can extend the $\zeta$ which is built above on a fundamental domain
 to $\mathcal T\setminus\mathcal H$ equivariantly. Since this extension coincides with the identity 
 on $\mathcal H$, it extends by $Id\vert_{Ox_1x_2}$. This is a continuous extension because any point 
 of the plane $Ox_1x_2$ adheres only to $\mathcal H\setminus Ox_1x_2$ when considering 
 the closure of $\mathcal T\setminus Ox_1x_2$. 

\end{demo}

\begin{prop} \label{F2} Let $X$ be a compact topological space, $0<{\tau_1}<{\tau_2}<1$ and $\hat \xi:X\times[0,{\tau_2}]\to X\times[0,1]$ be a topological embedding which is the identity on $X\times\{0\}$,  $X\times[0,\tau_1]\subset\hat \xi(X\times[0,{\tau_2}])$. Then

1. there is a homeomorphism $\hat\zeta:X\times[0,\tau_2]\to \xi(X\times[0,\tau_2])$ such that  $\hat \zeta$ is identity on $X\times[0,{\tau_1}]$ and is $\hat\xi$ on $X\times\{\tau_2\}$. 

2. if for a set ${B}\subset X$ the equality $\hat\xi|_{{ B}\times[0,\tau_2]}=id|_{{ B}\times[0,\tau_2]}$ is true then $\hat \zeta|_{{ B}\times[0,\tau_2]}=id|_{{ B}\times[0,\tau_2]}$.
\end{prop} 
\begin{demo} Let us choose $l\in(\tau_1,\tau_2)$ such that $X\times[0,l]\subset\hat \xi(X\times[0,{\tau_2}])$. Define a homeomorphism $\kappa:[\tau_1,1]\to[0,1]$ by the formula $$\kappa(t)=\cases{(x,\frac{l(t-\tau_1)}{l-\tau_1}),~t\in[\tau_1,l];\cr
(x,t),~t\in[l,1].\cr}$$ Let $\mathcal K(x,t)=(x,\kappa(t))$ on $X\times[\tau_1,1]$. Then the required homeomorphism can be defined by the formula $$\hat \zeta(x,t)=\cases{(x,t),~t\in[0,\tau_1];\cr
\mathcal K^{-1}\xi(\mathcal K((x,s)))),~s\in[\tau_1,\tau_2].\cr}$$ Property 2 automatically follows from this formula.
\end{demo}
\medskip

{We now collect some facts of geometric topology  in dimension 2 and  3 on which the proof 
of Proposition \ref{F1} is based}. We begin with the Sch\"onflies  Theorem (see Theorem 10.4 in \cite{moise}).

\begin{prop} \label{schoenflies}Every topological embedding of $\mathbb S^1$ into $\R^2$ is the restriction of a global homeomorphism of $\R^2$ which is the identity map outside some compact set of the plane.
\end{prop}

One can derive the Annulus Theorem in dimension 2{; we state and prove it in the only case 
which we use. The coordinates of $\R^2$ are $(x_1,x_2)$}. The unit closed disc in $\R^2$ is
 denoted by $ \D^2$; its boundary is $\SS^1$. 
The annulus $2\D^2\smallsetminus int(\D^2)$ is denoted by $\A$. Let $\II$ denote the arc 
$\{1\leq x_1\leq 2, x_2=0\}$.} 

\begin{lemm} \label{2-annulus} {Let $g: 2\SS^1\cup \II\to  \R^2\smallsetminus (0,0)$
 be a topological embedding which surrounds the origin
in the direct sense and has the next properties: $g(\II)\subset \{x_1>0\}$, the image $g(2\SS^1)$
avoids the circle $C$ of radius $\frac 32$ and $g(1,0)$ lies inside $\frac 32\D^2$. Then $g\vert _{2\D^2}$} extends to 
an embedding $G:\A\to \R^2\smallsetminus int(\D^2)$ which coincides with the identity on 
$\SS^1$  and maps $\II$ into $\{x_1>0\}$.
\end{lemm}
\begin{demo} Let $p$ be the last point on $g(\II)$ starting from $g(1,0)$ which  belongs to $C$.
 Let $q$ be its inverse image in $\II$. Define $G$ on the segment $[(1,0), q]$ as the affine map whose image is $[(1,0), p]$ and take $G$ coinciding with $g$ on $[q,(2,0)]$. The image $G(\II)$ is a simple arc in $\{x_1>0\}$. Because, any simple arc is tame in the plane, this definition of $G$ on $\II$
and the values which are imposed on the two circles $\SS$ and
 $2\SS^1$ extends to a neighborhood $N$ 
of $\SS^1\cup \II\cup 2\SS^1$ in $\R^2$. By taking one boundary component of $N$
one derives a parametrized simple curve $C'$ in $\R^2\smallsetminus \D^2$ which does not surround the origin. Therefore, by the Sch\"onflies Theorem, $C'$ bounds a disc $D$ in $\R^2\smallsetminus \D^2$
and the parametrization of $C'$ extends to a parametrization of $D$. This yields the complete
definition of $G$.
\end{demo}




We are now going to apply famous theorems of geometric topology in dimension 3 to
 a concrete situation emanating from 
the problem we are facing in Proposition \ref{F1}. The setting is the following.
We look at the 3-space 
$$Y= \A\times [0,1]=\{(x_1,x_2,x_3)\in\R^3\mid 1\leq x_1^2+x_2^2\leq 4,\, 0\leq x_3\leq 1\}.$$
Denote $Q_0$ the solid torus  in $Y$ limited by the next two annuli:

-- the vertical annulus $A_v:= \{ x_1^2+x_2^2= 1,\ 0\leq x_3\leq 1\}$;

-- the {\it standard curved} annulus $A_0:=  \{ x_1^2+x_2^2+ (x_3-\frac 12)^ 2= \frac 54\}\cap Y$.

\nd Observe that $A_0$ is contained in $Y$ and contains the two horizontal circles forming $\partial A_v$.

 \begin{prop}\label{solid} Let $g: A_0\to Y$ a \emph{bi-collared} (meaning that $g$ extends to 
 an embedding $A_0\times (-\ep, +\ep)$ for some $\ep>0$) topological embedding. It is assumed that
 $g$ coincides with the identity on $\partial A_0$ and maps the arc $\Ga_0:= A_0\cap \{x_1>0,x_2=0\}$
to an arc $\Ga$  in $Y\cap \{x_1>0\}$. Then, $g$ extends to an embedding $G: Q_0\to Y$
 which coincides with the identity on $A_v$.
 \end{prop}
 \begin{demo} 
 The image $A:=g(A_0)$ separates $Y$ in two components $X$ and $X^*$. Since $A$ is bi-collared,
 these two domains of $Y$ are topological 3-manifolds. Therefore, we are allowed to apply
 E. Moise's Theorem (See 
 \cite[Chap.23 \& Theorem 35.3]{moise} for existence and \cite[36.2]{moise} for uniqueness), including the {\it Hauptvermutung} in dimension 3:
 \vskip -.6cm
 ${}$
 \begin{quote}
{\it Every 3-manifold has a unique $PL$-structure, {up to an arbitrarily small topological isotopy}. 
 Moreover, its topological boundary is  a $PL$-submanifold.}
 \end{quote}
  As a consequence, $X$ and $X^*$ have $PL$-structures which agree on their intersection
  $A$. By uniqueness applied for $X\cup X^*$, the $PL$-structure on the union is the standard one
  after some  $C^0$-small ambient  isotopy.
  Denote by $P$ the planar surface  $\{x_2= 0, x_1>0\}\cap Y$. After a new $C^0$-small ambient isotopy
  in $Y$, we may assume that $A$ and $P$ are in general position. In what follows, we borrow
  the idea of proof from \cite[Theorem 3.1]{GZM}\footnote{In this article which we referred to, it should be
   meant that the $PL$-category (or smooth category)   is used. Indeed, there is no {\it general position} statement in topological geometry without more specific assumption.}.
  
  Since $P$ intersects each connected component of $\partial A$ in one point  only, we are sure that
  in general position $P\cap A$ is made of   finitely many simple closed curves $c_1,\ldots, c_k$ in $intA$ and one arc $\ga$ which links the two   components of $\partial A$. One of the above curves is {\it innermost} in $P$, meaning that   it bounds a disc in $P$ whose interior avoids $A$; let say that $c_1$ is so. More precisely, $c_1$ bounds a disc $d$ in $P$ and a disc $\de$ in $A$. By innermost position, $d\cup\de$
  is an embedded $PL$ 2-sphere $\si$. As $Y$ lies in $\R^3$, this sphere bounds a  3-ball 
  $\De\subset Y$. We are going to use these data in two ways.
  
 First, we use $\de$ for finding an isotopy $h_t$ of $A$ in itself  from $Id\vert _A$
 to $h_1: A\to A$ such that $h_1(\Ga)\cap \de= \emptyset$. This is easily done as $\Ga\cap\de$
 avoids one point $z_\de$ in $\de$: one pushes $\Ga\cap\de$ along the rays of $\de$ issued $z_\de$. 
 Notice that $h_1(\Ga)$ still lies in $\{x_1>0\}$, but this could be no longer true for $h_t(\Ga), \ t\neq 0,1,$
 when $\de$ is not contained in $\{x_1>0\}$.
 
 Once, this is done, the ball $\De$ is used for finding an ambient isotopy of $Y$ which is supported in
 a neighborhood of $B$, small enough so that $h_1(\Ga)$ is kept fixed, and which moves
 $A\cap \De$ to the complement of $P$. Hence, this isotopy cancels $c_1$ from $A\cap P$; all
 intersection curves contained in $int\,\de$ are cancelled at the same time. By repeating
  isotopies similar to the two previous ones, as many times as necessary, we get an embedding
  $g': A_0\to Y$ which coincides with the identity on $\partial A_0$, still maps 
  $\Ga_0$ into $\{x_1>0\}$ and fulfills the following property:
  $A':= g'(A_0)\cap P$ is made of one arc $\ga'$ only which links the two components 
  of $g'(\partial A_0)$. The annulus $A'$ divides $Y$ into two (closed) domains $X'$ and $X'^*$ which 
  come from the splitting $X\cup_A X^*$ by an ambient isotopy fixing $\{x_3=0,1\}$ pointwise.
  
  As $\ga'$ is the only intersection component of 
  $P\cap A'$, one knows  that $\ga'$ divides $P$ into a disc $\mu'\subset X'$ 
  (meaning a {\it meridian} disc in a solid torus) and its 
  complement in $P$. Removing from
   $X'$ a regular neighborhood of
  $\mu'$ yields a $PL$ embedded 2-sphere $S$. 
  According to the Alexander theorem \cite{alexander},
  this sphere bounds a ball  $B_{X'}$ in $\R^3$, as $Y\subset \R^3$. 
   It is not possible that $B_{X'}$ contains $\mu'$ in its interior; in the contrary, $B_{X'}$ 
   would get out of $X'^*$ and have a non-bounded interior. 
   As a consequence, $X'$ is a \emph{solid torus} since it is made of a ball and a 1-handle attached.
   The same holds for $X$ as it is ambient isotopic to $X'$ in $Y$. 
   
   This is not sufficient 
   for concluding. It would be necessary to prove the same for the curve $g'(\Ga_0)$,
   after  making it  a closed curve 
   by adding the vertical arc $\ga_0^*\subset A_v$ which links the two points 
   of $g'(\Ga_0)\cap A_v=
   \Ga_0\cap A_v= g(\Ga_0)\cap A_v= \Ga_0\cap A_v= \ga'\cap A_v$.
   That is the place where the assumption about $g(\Ga_0)$ is used.\\

  \nd{\sc Claim.} {\it There exists an ambient isotopy  
  from $Id\vert _Y$ to $k_1$
  which is stationary on the vertical annulus $A_v$, 
  which maps $A'$ into itself and moves  $g'(\Ga_0)$ to $\ga'$.}\\
  
  \nd{\sc Proof of the claim.} Assume first that the arcs $g'(\Ga_0)$ and $\ga'$ meet in their end points only.
   Consider the closed curve $\al$ which is  made of  $\ga\cup g'(\Ga_0)$; it is
   contained in $\{x_1>0\}\cap A'$. 
  By construction, the homological intersection of $\al$ with $\ga'$ is zero. Therefore, $\al$ bounds a disc
  $\de'\subset A'$. Notice that it could be not contained in $\{x_1>0\}$. The disc $\de'$
  allows one to move $g'(\Ga_0)$ to $\ga'$ by an isotopy of $A'$ into itself
  with the required properties. 
  
  In case where $g'(\Ga_0)\cap int\ga'$ is non-empty, in general position this intersection is made of finitely 
  many points. Among them, choose the point $x$ which is the closest to $\ga'\cap \{x_3=1\}$ when 
  traversing $\ga'$ starting from  bottom. Denote by $x_0$ the point $\ga'\cap \{x_3= 1\}$.
  One forms a closed curve $\tau\subset A'$ made of two arcs in $A'$,
  from $x$ to $x_0$ respectively in $\ga'$ and $g'(\Ga_0)$.
 For the same reason as for $\al$ above,
 the curve $\tau$ bounds a disc in $A'$  which allows one to cancel $x$  from 
 $g'(\Ga_0)\cap\ga'$ by an isotopy of $A'$ into itself with the required properties.
 Iterating this process reduces us to the first case. The claim is proved. $\diamond$\\

 Let $g'': = k_1 g': A_0\to A'$. As a consequence of the claim, the closed curve $g''(\Ga_0)\cup \ga_0^*$
 is the boundary of the {\it meridian} disc $\mu'$.
 We are going to show that  $g''$ extends to an embedding
 $G'': Q_0\to X'$ which coincides with the identity on $A_v$. In this aim, we denote by 
 $\mu_0$ the meridian of $Q_0$ defined by $\mu_0= Q_0\cap P$. For beginning with,
  we consider regular neighborhoods $N(\mu_0)$ and $N(\mu')$ of both meridians and we extend 
 $g''\vert_{N(\mu_0)\cap A_0}$ to a homeomorphism
  $G''_0: N( \mu_0)\to N(\mu')$ which is the identity on $N( \mu_0)\cap A_v$. 
   
   Let $B_{Q_0}$ be the ball in $Q_0$ which is the closure of $Q_0\setminus N(\mu_0)$.
   The restricted map $G''_0\vert_{\partial B_{Q_0}}$ glued with 
   the restriction of $g''$ to the closure of $A_0\setminus N(\mu_0)$ yields a homeomorphism
   $$G''_1: \partial B_{Q_0}\to \partial  B_{X'}.$$ 
   The desired $G''$ is obtained by extending
   $G''_1$ to $B_{Q_0}$ by the cone construction (seeing a ball as the cone on its boundary).
   Since $g$ and $g''$ are related one to the other by an ambient isotopy fixing $A_v$ pointwise,
   an extension of $g$ follows from an extension of $g''$.
   \end{demo}

{\nd address: {Laboratoire de Topologie, UMR 5584 du CNRS, Dijon, France}

\nd email: {bonatti@u-bourgogne.fr}}
\vspace{10mm}

{\nd address: {National Research University Higher School of Economics,  603005, Russia, Nizhny Novgorod, B. Pecherskaya, 25}

\nd email: {vgrines@yandex.ru}}
\vspace{10mm}

{\nd address: {Universit\'e de Nantes, LMJL, UMR 6629 du CNRS, 44322 Nantes, France}

\nd email: {francois.laudenbach@univ-nantes.fr}}
\vspace{10mm}

{\nd address: {National Research University Higher School of Economics,  603005, Russia, Nizhny Novgorod, B. Pecherskaya, 25}

\nd email: {olga-pochinka@yandex.ru}}

\end{document}